\documentclass[a4paper,12pt]{amsart}
\usepackage{amsmath,amsthm,amssymb}
\usepackage{amsfonts}
\usepackage{txfonts}
\usepackage{latexsym}
\usepackage{iftex}
\ifpTeX
\usepackage[dvipdfmx]{graphics,color}
\else
\usepackage{graphicx,color}
\fi
\usepackage{eucal}
\usepackage{mathrsfs}
\ifpTeX
\usepackage[dvipdfmx]{hyperref}
\usepackage{pxjahyper}
\else
\usepackage{hyperref}
\fi
\usepackage{tikz-cd}

\newcommand{\Gm}{\mathbb{G}_m}

\theoremstyle{plain}
\newtheorem{theorem}{Theorem}[section]
\newtheorem{proposition}[theorem]{Proposition}
\newtheorem{lemma}[theorem]{Lemma}
\newtheorem{corollary}[theorem]{Corollary}

\theoremstyle{definition}
\newtheorem{definition}[theorem]{Definition}

\makeatletter
\renewenvironment{proof}[1][\proofname]{\par
  \normalfont
  \topsep6\p@\@plus6\p@ \trivlist
  \item[\hskip\labelsep{\bfseries #1}\@addpunct{\bfseries.}]\ignorespaces
}{%
  \endtrivlist
}
\renewcommand{\proofname}{proof}
\theoremstyle{remark}
\newtheorem{remark}[theorem]{Remark}

\numberwithin{equation}{section}

\title{Algebraic cobordism via spans}
\date{\today} 
\author{Yuki Kato}
\thanks{The author is supported by Grants-in-Aid for Scientific Research No. 23K03080, Japan Society for the Promotion of Science.}

\address{National Institute of Technology, Kurume College, 
	      1-1-1, Komorono, Kurume, Fukuoka, 830-8555 JAPAN.}
\email{kato\_051@kurume-nct.ac.jp}
\subjclass[2020]{14F42 (primary), 14G45, 18N40}
\keywords{algebraic cobordism, spans, projective bundle formula, perfectoid tilting}

\newcommand{\MGL}{\mathrm{MGL}}
\newcommand{\Map}{\mathrm{Map}}

\newcommand{\Sp}{\mathrm{Sp}}             
\newcommand{\Span}{\mathbf{Span}}          

\usepackage{enumerate}
\usepackage{array}
\usepackage[all]{xy} 
\usepackage{delarray}

\begin{document}
\begin{abstract}
We define the algebraic cobordism of $\infty$-categories equipped with universal line bundle data as an initial oriented functor in the associated span category. In the standard motivic framework, this recovers the Thom spectrum model established by Voevodsky; Gepner and Snaith. Furthermore, assuming that the $\infty$-category contains Grassmann objects of all ranks, we prove that the projective bundle formula and the corresponding Chern-class and Whitney-sum identities hold for any oriented functor satisfying the splitting principle property.

We apply the span formalism to perfectoid geometry. For perfectoid algebras $R$ with tilt $R^\flat$, we construct perfectoid cobordism, prove tilting equivalences, and compare the arc-local and $v$-local $p$-adic theories.
\end{abstract}
\maketitle

\theoremstyle{plain}
\newtheorem{maintheorem}[theorem]{Main Theorem}
\newtheorem{assumption}[theorem]{Assumption}
\section{Introduction}
\label{sec:introduction}
Voevodsky~\cite{V} introduced the $\mathbb{P}^1$-stable homotopy category
$\mathrm{SH}(S)$ for any scheme $S$ as an algebraic analogue of classical stable
homotopy theory. In his framework, the algebraic cobordism $\MGL$ serves as the universal oriented motivic spectrum, which is analogous to the complex cobordism $\mathrm{MU}$.  Panin~\cite{Panin2003Oriented} constructed orientations for cohomology theories using correspondences that rely on span-theoretic functoriality. Following Voevodsky's approach, Garkusha and Panin~\cite{GP} developed a theory of framed motives based on framed correspondences, which consist of spans equipped with finite syntomic trivializations of their normal bundles. Subsequently, Elmanto, Hoyois, Khan, Sosnilo, and Yakerson~\cite{ModMGL} constructed a stable $\infty$-category whose monoidal unit is $\MGL$ by introducing finite syntomic correspondences of derived schemes. 

In this paper, we formalize these geometric constructions using the $\infty$-categorical framework of spans explained by Hoyois~\cite{Hoyois2017SixOps}. This perspective enables us to construct Gysin morphisms $\infty$-categorically.

This paper consists of two parts. The first part constructs the algebraic
cobordism $\MGL_\mathcal{C}$ for an $\infty$-category
$(\mathcal{C},\mathbf{S})$ equipped with universal line bundle data, 
thereby defining $\MGL_\mathcal{C}$ as the universal oriented spectrum for the
span category $\Span_{\mathbf{S}}(\mathcal{C})$. Under the
hypotheses of Section~\ref{sec:pbf}, we prove the projective bundle
formula and the identity for the Whitney sums. 
In the context of Morel and Voevodsky's~\cite{MV}, this recovers Voevodsky's algebraic cobordism spectrum $\MGL$~\cite{V}.

The second part applies the span formalism to Scholze's~\cite{Scholze2012Perfectoid} perfectoid geometry. We work with
$\mathrm{Stk}_v(\mathrm{Perf}_R)$, where a $v$-stack means a $v$-sheaf on the
$\infty$-category of perfectoid spaces for perfectoid algebras $R$. In this
setting, quasi-smooth morphisms are unsuitable, as 
the $p$-adic complete cotangent complexes are acyclic for perfectoid rings. In
Section~\ref{sec:perfectoid}, we define the admissible class of morphisms as being
generated by smooth morphisms and zero-sections of vector bundles.

We present the main results. 

\begin{maintheorem}[Theorem~\ref{thm:pbf-universal} and Corollary~\ref{cor:pbf-general}]\label{thm:intro-A}
Let $E\colon \Span_{\mathbf{S}}(\mathcal{C})\to\Sp$ be an oriented
cohomology theory on an $\infty$-category $(\mathcal{C},\mathbf{S})$
equipped with universal line bundle data. Assume that the underlying
$\infty$-category is locally presentable, that $\mathbf{BG}_m$ is
group-like, and that the induced morphism
\[
\Phi \colon \mathfrak{L}^\circ \times \mathbb{A}^1 \to \mathfrak{L}^\circ \times_{\mathbf{BG}_m} \mathfrak{L}
\]
is a weak equivalence whose restriction $\Phi^\circ$ makes
$\mathfrak{L}^\circ\to\mathbf{BG}_m$ a principal $\mathbb{G}_m$-torsor and
whose base changes used in Section~\ref{sec:pbf} remain principal
$\mathbb{G}_m$-torsors. If $V\to X$ is a split rank-$n$ vector bundle in
$\mathcal{C}$, then the canonical map
\[
\Phi_V\colon \bigoplus_{i=0}^{n-1}E(X)\cdot\xi^i\longrightarrow E(\mathbb{P}(V))
\]
is a weak equivalence of $E(X)$-modules. If, in addition, $\mathcal{C}$ has Grassmann objects of all ranks and $E$ satisfies the splitting principle property of
Definition~\ref{prop:bgl-assumption}, then for any vector bundle
$V \to X$, the same canonical map $\Phi_V$ is a weak
equivalence.
\end{maintheorem}

In the setting of perfectoid geometry, the geometric input is
given by Definition~\ref{def:admissible-class}, which specifies the admissible
class.
\begin{maintheorem}[Theorem~\ref{thm:tilting-mglperf}]\label{thm:intro-B}
For any perfectoid algebra $R$ with tilt $R^\flat$, there is a
canonical equivalence
\[
\MGL_{\mathrm{perf}}(R)\simeq\MGL_{\mathrm{perf}}(R^\flat),
\]
of $\mathbb{E}_\infty$-ring spectra, which is functorial for perfectoid algebras over $R$.
\end{maintheorem}

For the $p$-adic theory, write
\[
\MGL_\mathrm{perf}^{v}(-,\mathbb{Z}_p)=L_v\bigl(\MGL_{\mathrm{perf}}(-,\mathbb{Z}_p)\bigr),
\qquad
\MGL_\mathrm{perf}^{\mathrm{arc}}(-,\mathbb{Z}_p)=L_{\mathrm{arc}}\bigl(\MGL_{\mathrm{perf}}(-,\mathbb{Z}_p)\bigr)
\]
for the $v$-sheafification and arc-sheafification.

\begin{maintheorem}[Corollary~\ref{thm:tilting-zp} and Theorem~\ref{thm:arc-v-comparison}]\label{thm:intro-C}
Let $\mathrm{Perf}$ denote the category of perfectoid spaces, and let $\MGL_{\mathrm{perf}}^{v}(-, \mathbb{Z}_p)$ be the $p$-adic complete algebraic cobordism presheaf on $\mathrm{Perf}$. Then the following hold:
\begin{enumerate}
    \item For any perfectoid algebra $R$ with tilt $R^\flat$, there is a canonical equivalence of $\mathbb{E}_\infty$-ring spectra:
    \[
    \MGL_{\mathrm{perf}}^{v}(R, \mathbb{Z}_p) \simeq \MGL_{\mathrm{perf}}^{v}(R^\flat, \mathbb{Z}_p),
    \]
    which is strictly functorial with respect to perfectoid $R$-algebras.
    \item For any perfectoid algebra $R$, the natural comparison morphism induces a canonical equivalence:
    \[
    \MGL_{\mathrm{perf}}^{\mathrm{arc}}(R, \mathbb{Z}_p) \xrightarrow{\sim} \MGL_{\mathrm{perf}}^{v}(R, \mathbb{Z}_p),
    \]
    implying that $\MGL_{\mathrm{perf}}^{v}(-, \mathbb{Z}_p)$ satisfies hypercomplete $\mathrm{arc}$-descent on perfectoid rings.
\end{enumerate}
\end{maintheorem}

The paper is organized as follows: Section~\ref{sec:spans} introduces universal line bundle data,
admissible classes of morphisms, and the span formalism.
Section~\ref{sec:mgl-construction} defines $\MGL_\mathcal{C}$ as the
universal oriented spectrum. Section~\ref{sec:pbf} establishes the projective
bundle formula and the identities for Chern classes and Whitney sums.
In Section~\ref{sec:perfectoid}, we introduce the algebraic cobordism $\MGL_{\mathrm{perf}}$ of perfectoid algebras and prove the tilting equivalences.

\section{\texorpdfstring{The $\infty$-category of spans and orientations}{The infinity-category of spans and orientations}}
\label{sec:spans}

This section introduces $\infty$-categories with universal line bundle data, explores their $\mathbb{A}^1$-homotopy theory, and discusses span categories that include admissible classes of morphisms and oriented cohomology functors arising from these constructions.
   
\subsection{Universal line bundle data and orientation}
\begin{definition}
Let $\mathcal{C}$ be an $\infty$-category that admits finite limits and
colimits. Universal line bundle data consist of a triple
$(\pi:\mathfrak{L}  \to\mathbf{BG}_m,\, i: \mathfrak{L}^\circ \to \mathfrak{L}, z:\mathbf{BG}_m \to \mathfrak{L})$ of morphisms
as follows:
\begin{itemize}
\item A morphism $\pi: \mathfrak{L} \to\mathbf{BG}_m$ in $\mathcal{C}$, referred to as the {\it universal line bundle object}.
\item A morphism $z:\mathbf{BG}_m \to \mathfrak{L}$, which satisfies the condition that $\pi \circ z$ is homotopic to the identity, is called the {\it zero section} of $\pi$.
\item A morphism $i: \mathfrak{L}^\circ \to \mathfrak{L}$, called the {\it punctured line object}, such that the pull-back of $i$ along $z$ is weakly equivalent to the initial object $\emptyset$.
\item The object $\mathbf{BG}_m$ is a commutative monoid object of $\mathcal{C}$ with respect to Cartesian products. 
Specifically, there exists a multiplication map:
\[
\mu:\mathbf{BG}_m \times\mathbf{BG}_m \to\mathbf{BG}_m. 
\]
\end{itemize}
\end{definition}
For any morphism $f: X \to\mathbf{BG}_m$, we define the line bundle, the punctured line bundle, and the {\it Thom space} as follows:
\begin{itemize}
\item The {\it line bundle} is $\mathcal{L}_X = X \times_{\mathbf{BG}_m} \mathfrak{L} \to X$. 
\item The {\it punctured line bundle} is $\mathcal{L}_X^\circ = X \times_{\mathbf{BG}_m} \mathfrak{L}^\circ \to X$. 
\item The {\it Thom space} $\mathrm{Th}(\mathcal{L}_X)$ is the homotopy cofiber $\mathrm{Th}(\mathcal{L}_X) = \mathrm{Cofib}(\mathcal{L}_X^\circ \to \mathcal{L}_X)$. 
\end{itemize}

By \cite[Proposition 2.4.3]{WY2024}, the cofiber functor $\mathrm{Cofib}: \mathrm{Ar}(\mathcal{C}) \to \mathcal{C}$ is a symmetric monoidal functor, where the arrow category $\mathrm{Ar}(\mathcal{C})$ is equipped with the pushout product 
\[
    (f:X_0 \to X_1)\Box (g:Y_0 \to Y_1) = (X_0 \times Y_1)\amalg_{X_0 \times Y_0} (X_1 \times Y_0) \to X_1 \times Y_1 
    \]
as its monoidal structure.
\begin{lemma}
\label{lem:thom-symmetric-monoidal}
The Thom-space functor $\mathrm{Th}:\mathcal{C}_{/\mathbf{BG}_m} \to \mathcal{C}$ is symmetric monoidal, compatible with the pushout product. More precisely, for line objects $\mathcal{L}_X, \mathcal{L}_Y \in \mathcal{C}_{/\mathbf{BG}_m}$, one has a canonical equivalence in $\mathcal{C}$:
\[
\mathrm{Th}(\mathcal{L}_X\,\Box\ \mathcal{L}_Y)\simeq \mathrm{Th}(\mathcal{L}_X)\wedge \mathrm{Th}(\mathcal{L}_Y),
\]
where $\mathcal{L}_X\,\Box\ \mathcal{L}_Y$ denotes the pushout product of the respective punctured unit disk inclusions $i_X\colon \mathcal{L}_X^\circ \to \mathcal{L}_X$ and $i_Y\colon \mathcal{L}_Y^\circ \to \mathcal{L}_Y$:
\[
i_X\,\Box\ i_Y = (\mathcal{L}_X^\circ \times \mathcal{L}_Y) \amalg_{\mathcal{L}_X^\circ \times \mathcal{L}_Y^\circ } (\mathcal{L}_X \times \mathcal{L}_Y^\circ) \to \mathcal{L}_X \times \mathcal{L}_Y.
\]
\end{lemma}

\begin{proof}
By definition, the Thom space $\mathrm{Th}(\mathcal{L})$ is the cofiber of the inclusion $i\colon \mathcal{L}^\circ \to \mathcal{L}$. Since the assignment $\mathcal{L} \mapsto (\mathcal{L}^\circ \to \mathcal{L})$ defines a symmetric monoidal functor from $\mathcal{C}_{/\mathbf{BG}_m}$ to the symmetric monoidal $\infty$-category $(\mathrm{Ar}(\mathcal{C}), \Box)$, and the cofiber functor $\mathrm{Cofib}: \mathrm{Ar}(\mathcal{C}) \to \mathcal{C}$ is symmetric monoidal by \cite[Proposition 2.4.3]{WY2024}, their composition $\mathrm{Th}$ is also symmetric monoidal. \qed
\end{proof}

Given universal line bundle data and a rational point $x: * \to\mathbf{BG}_m$, we define the {\it affine line} $\mathbb{A}^1$, the {\it algebraic circle}
$\mathbb{G}_m$, and the {\it projective line} $\mathbb{P}^1$ as follows:
\[
 \mathbb{A}^1 = x^*(\mathfrak{L}), \quad \mathbb{G}_m = x^*(\mathfrak{L}^\circ),  \quad \mathbb{P}^1 = \mathrm{Th}(\mathbb{A}^1) 
\]
in the $\infty$-category $\mathcal{C}$. These objects define the $\mathbb{A}^1$-homotopy theory of
$\mathcal{C}$.
\begin{definition}
Let $\mathcal{C}$ be an $\infty$-category equipped with an affine line object $\mathbb{A}^1$.
\begin{enumerate}
\item An object $Z \in \mathcal{C}$ is said to be {\it $\mathbb{A}^1$-local} if, for any object $X \in \mathcal{C}$, 
the projection $p: X \times \mathbb{A}^1 \to X$ induces an equivalence:
\[
p^*: \mathrm{Map}_{\mathcal{C}}(X, Z) \to \mathrm{Map}_{\mathcal{C}}(X \times \mathbb{A}^1, Z).
\]
\item A morphism $f: X \to Y$ in $\mathcal{C}$ is called an {\it $\mathbb{A}^1$-weak equivalence} 
if for any $\mathbb{A}^1$-local object $Z$, the induced map
\[
f^*: \mathrm{Map}_{\mathcal{C}}(Y, Z) \to \mathrm{Map}_{\mathcal{C}}(X, Z)
\]
is a homotopy equivalence.
\end{enumerate}
\end{definition}

\subsection{Spans with the admissible class of morphisms}

To define the $(\infty,\,2)$-category of spans, we first specify the admissible class of morphisms.
\begin{definition}\label{def:axiomatic-qsm-class}
Let $\mathcal{C}$ be an $\infty$-category admitting finite limits and
colimits with universal line bundle data $(\pi:\mathfrak{L} \to\mathbf{BG}_m,\, i: \mathfrak{L}^\circ \to
\mathfrak{L}, z:\mathbf{BG}_m \to \mathfrak{L})$ and a rational point
$x: * \to\mathbf{BG}_m$. Let $\mathbf{S}$ be a class of morphisms in $\mathcal{C}$.  We call it the class of {\it admissible morphisms}  
if it satisfies the following conditions:
\begin{enumerate}
\item The class $\mathbf{S}$ contains all weak equivalences, is closed under composition, and is stable under base changes.
\item For any object $X \in \mathcal{C}$, the codiagonal $\nabla: X \amalg X \to X$ belongs to $\mathbf{S}$.
\item The morphisms $\pi$ and $z$ belong to $\mathbf{S}$, and $\pi$ is an $\mathbb{A}^1$-equivalence.
\end{enumerate}
\end{definition}
\begin{definition}
Let $\mathcal{C}$ be an $\infty$-category equipped with an admissible class
of morphisms $\mathbf{S}$. The $(\infty,\,2)$-category
$\mathrm{Span}_{\mathbf{S}}(\mathcal{C})$, called the
\textit{$\mathbf{S}$-span of $\mathcal{C}$}, has the following:
\begin{itemize}
    \item The objects of $\mathrm{Span}_{\mathbf{S}}(\mathcal{C})$ are the
    same as those of $\mathcal{C}$.
    \item A 1-morphism from $X$ to $Y$ is a diagram
    $[X \xleftarrow{s} Z \to Y]$, where $s$ belongs to $\mathbf{S}$. The
    composition of $[X \xleftarrow{s} V \to Y]$ and
    $[Y \xleftarrow{t} W \to Z]$ is defined by the homotopy Cartesian
    product, giving the diagram:
    \[
    \xymatrix{
      & & V \times_Y W \ar[dl]_{s'} \ar[dr]^{t'} & & \\
      & V \ar[dl]_{s} \ar[dr] & & W \ar[dl]^{t} \ar[dr] & \\
    X & & Y & & Z.
    }
    \]
    (Since $\mathbf{S}$ is an admissible class, it is stable under
    pull-backs and composition, so $s \circ s' \in \mathbf{S}$.)
    \item A 2-morphism from $[X \leftarrow Z \to Y]$ to
    $[X \leftarrow W \to Y]$ is a morphism $\alpha \colon Z \to W$ in
    $\mathcal{C}$ that makes the following diagram commute up to homotopy:
    \[
    \xymatrix{
    & Z \ar[dl] \ar[dr] \ar[dd]^{\alpha} & \\
    X & & Y \\
    & W \ar[ul] \ar[ur] &
    }
    \]
\end{itemize}
\end{definition}

\begin{definition}
For any morphism $f: X \to Y$ in $\mathbf{S}$, let $f^*: Y \to X$ denote
the induced span
$[Y \xleftarrow{f} X \xrightarrow{\mathrm{id}_X} X ]$, called the
{\it pull-back} of $f$, and let $f_!$ denote
$[X \xleftarrow{\mathrm{id}_X} X \xrightarrow{f} Y]$, called the
{\it (Gysin) push-forward} of $f$.

The 2-morphisms in $\mathrm{Span}_{\mathbf{S}}(\mathcal{C})$ given by the
following commutative diagrams
\[
\xymatrix@1{
& X \ar[dl]_{\mathrm{id}_X} \ar[dr]^{\mathrm{id}_X} \ar[dd]^{\mathrm{diag}}  & \\
X & & X \\
& X \times_Y X \ar[ul]^{\mathrm{pr}_1} \ar[ur]_{\mathrm{pr}_2} &,
}
\]
and
\[
\xymatrix@1{
& X \ar[dl]_f \ar[dr]^f \ar[dd]^f  & \\
Y & & Y \\
& Y \ar[ul]^{\mathrm{id}_Y} \ar[ur]_{\mathrm{id}_Y} &
}
\]
determine, respectively, the unit
$\eta \colon \mathrm{id}_X \Rightarrow f^*\circ f_!$ and the counit
$\epsilon \colon f_! \circ f^* \Rightarrow \mathrm{id}_Y$.
\end{definition}
If the admissible class $\mathbf{S}$ is stable under finite coproducts, then the span category $\Span_{\mathbf{S}}(\mathcal{C})$ has the canonical additive structure defined by the codiagonals.
\begin{definition}
\label{prop:span_addition}
Let $\mathcal{C}$ be an $\infty$-category with finite coproducts, and
assume that the admissible class $\mathbf{S}$ is stable under coproducts. Then
the hom-spaces in the $(\infty, 2)$-category
$\mathrm{Span}_{\mathbf{S}}(\mathcal{C})$ naturally carry a commutative
monoid structure, called the {\it addition}. For any
$\alpha = [A \xleftarrow{s_1} Z_1 \xrightarrow{t_1} X]$ and
$\beta = [A \xleftarrow{s_2} Z_2 \xrightarrow{t_2} X]$, $\alpha+\beta$ is
defined by
\[
    \alpha + \beta \quad = \quad (\nabla_X)_! \circ (\alpha \amalg \beta) \circ (\nabla_A)^*.
\]
Associativity and commutativity of the addition follow from the
universal properties of coproducts in $\mathcal{C}$.
\end{definition}

The $(\infty,\,2)$-category $\Span_\mathbf{S}(\mathcal{C})$ has the projection formula property with respect to the following intersection product.  
\begin{definition}
\label{def:tensor_product}
Let $\mathcal{C}$ be an $\infty$-category with finite limits. For any
objects $A, B \in \mathcal{C}_{/X}$, their exterior product $A \times B$ is
naturally an object over $X \times X$. We define the relative tensor
product $A \cdot_X B$ geometrically in
$\mathrm{Span}_{\mathbf{S}}(\mathcal{C})$ as the pull-back of the exterior
product along the diagonal map $\Delta_X \colon X \to X \times X$:
\begin{equation}
    A \cdot_X B \quad = \quad \Delta_X^* (A \times B)
\end{equation}
where the pull-back operator is the span
$\Delta_X^* = [X \times X \xleftarrow{\Delta_X} X \xrightarrow{\mathrm{id}_X} X]$.
\end{definition}
 
\begin{proposition}
\label{prop:span_projection}
Let $\mathcal{C}$ be an $\infty$-category with finite limits and an admissible class $\mathbf{S}$. For any morphism $f \colon X \to Y$ in $\mathbf{S}$ and spans
$\alpha = [X \xleftarrow{a} A \xrightarrow{p} Z_A]$ over $X$ and
$\beta = [Y \xleftarrow{b} B \xrightarrow{q} Z_B]$ over $Y$, there is a canonical
equivalence in
$\mathrm{Span}_{\mathbf{S}}(\mathcal{C})$:
\[
    f_!(\alpha \cdot_X f^*(\beta)) \simeq f_!(\alpha) \cdot_Y \beta.
\]
\end{proposition}

\begin{proof}
First, the pull-back $f^*(\beta)$ is represented by a span
over $X$:
\[
f^*(\beta) = [X \xleftarrow{\mathrm{pr}_1} X \times_Y B \xrightarrow{q \circ \mathrm{pr}_2} Z_B]
\]
The roof of $\alpha \cdot_X f^*(\beta)$ is the homotopy fiber product
$A \times_X (X \times_Y B)$, which is canonically equivalent to $A \times_Y B$. Therefore, we have an equivalence 
\[
\alpha \cdot_X f^*(\beta) \simeq [X \xleftarrow{a \circ \pi_A} A \times_Y B \xrightarrow{p \times q} Z_A \times Z_B]
\]
Applying the push-forward $f_!$, the left-hand side span over $Y$ is represented by
\[
f_!(\alpha \cdot_X f^*(\beta)) \simeq [Y \xleftarrow{f \circ a \circ \pi_A} A \times_Y B \xrightarrow{p \times q} Z_A \times Z_B]
\]

We compute the right-hand side $f_!(\alpha) \cdot_Y \beta$, where the push-forward is represented by a span
\[
f_!(\alpha) = [Y \xleftarrow{f \circ a} A \xrightarrow{p} Z_A].
\]
The roof is the homotopy fiber product $A \times_Y B$ of $A \xrightarrow{f \circ a} Y$ and $B \xrightarrow{b} Y$, with projections $\pi_A$ and $\pi_B$. Since $(f \circ a) \circ \pi_A = b \circ \pi_B$, this gives
\[
f_!(\alpha) \cdot_Y \beta \simeq [Y \xleftarrow{f \circ a \circ \pi_A} A \times_Y B \xrightarrow{p \times q} Z_A \times Z_B].
\]
Hence, we obtain that $f_!(\alpha \cdot_X f^*(\beta)) \simeq f_!(\alpha) \cdot_Y \beta$. \qed
\end{proof}

\begin{theorem}[Projection Formula]
\label{thm:universal_projection}
Let $\mathcal{C}$ be an $\infty$-category with finite limits, and
let $\mathbf{S}$ be an admissible class.
Given any symmetric monoidal functor
$E \colon \mathrm{Span}_{\mathbf{S}}(\mathcal{C}) \to \Sp$ and any morphism $f:X\to Y$ in $\mathbf{S}$, there is a canonical equivalence
\begin{equation}
    f_!(a \cdot f^*(b)) \simeq f_!(a) \cdot b
\end{equation}
for any classes $a \in E(X)$ and $b \in E(Y)$.
\end{theorem}
\begin{proof}
This follows from Proposition~\ref{prop:span_projection} and the $\infty$-categorical Yoneda lemma~\cite[p.~317, Proposition 5.1.3.1]{HT}. \qed
\end{proof}

\section{Oriented spectrum and the algebraic cobordism of \texorpdfstring{$\infty$}{infinity}-categories with universal line bundle data}
\label{sec:mgl-construction}
In this section, we define the algebraic cobordism $\MGL_\mathcal{C}$ of the $\infty$-category $\mathcal{C}$, equipped with universal line-bundle data, as the universal oriented spectrum.

\subsection{\texorpdfstring{Oriented cohomology theory of $\infty$-categories}{Oriented cohomology theory of infinity-categories}}
By the definition of universal line bundle data from the previous section, the object $\mathbf{BG}_m$ is a commutative monoid object in $\mathcal{C}$. For any symmetric monoidal functor $E:\Span_\mathbf{S}(\mathcal{C}) \to \Sp$, the spectrum $E(\mathbf{BG}_m)$ has an $\mathbb{E}_\infty$-ring structure. 
We introduce oriented cohomology theories on the $(\infty,\,2)$-category of $\mathbf{S}$-spans.
\begin{definition}
\label{def:oriented-cohomology}
Let $\mathcal{C}$ be an $\infty$-category equipped with universal
line bundle data, and let $\mathbf{S}$ be an admissible class of
morphisms in $\mathcal{C}$.
A symmetric monoidal functor $E: \mathrm{Span}_{\mathbf{S}}(\mathcal{C}) \to \mathrm{Sp}$ is called an {\it oriented cohomology theory} if it satisfies the following
conditions:
\begin{enumerate}
\item For any object $X$, the projection $X \times \mathbb{A}^1 \to X$ induces an equivalence $E(X) \to E(X \times \mathbb{A}^1)$.
\item For any pair of objects $(X,\,Y)$, the canonical map induces an equivalence:
\[
E(X \amalg Y) \xrightarrow{\sim} E(X) \oplus E(Y).
\]
\item The underlying presheaf $E|_{\mathcal{C}^{\mathrm{op}}}: \mathcal{C}^{\mathrm{op}} \to \mathrm{Sp}$ is excisive, sending homotopy coCartisan squares in $\mathcal{C}$ to homotopy Cartesian squares of spectra.
\item Let $p:\mathfrak{L} \to \mathbf{BG}_m$ be the universal line bundle, and $z: \mathbf{BG}_m \to \mathfrak{L}$ its zero-section.  The Gysin push-forward $z_!$ induces an equivalence of $E(\mathbf{BG}_m)$-modules:
    \[
        E(\mathbf{BG}_m) \xrightarrow{z_!} E(\mathrm{Th}(\mathfrak{L})).
    \]
\end{enumerate}
\end{definition}

\begin{definition}
\label{def:universal-chern}
Let $\mathcal{C}$ be an $\infty$-category equipped with universal line bundle data, including the projection $\pi\colon \mathfrak{L} \to \mathbf{BG}_m$, its zero-section $z \colon \mathbf{BG}_m \to \mathfrak{L}$, and an admissible class of morphisms $\mathbf{S}$. 
Let $E$ be an oriented cohomology theory on $\mathcal{C}$. Since $z$ is admissible, it induces a Gysin push-forward $z_! \colon E(\mathbf{BG}_m) \to E(\mathfrak{L})$. The canonical \emph{Thom class} of the universal line bundle is defined as the push-forward of the unit element:
\[
    \mathrm{th}_\mathfrak{L} = z_!(1) \in E(\mathfrak{L}).
\]
Let $c^E_1$ denote the \emph{first Chern class} of the universal line bundle $\mathfrak{L}$, defined as the restriction of the Thom class along the zero-section $z$:
\[
    c_1^E = z^*(\mathrm{th}_\mathfrak{L}) = z^* \circ z_! (1) \in E(\mathbf{BG}_m).
\]
(Note that the induced map $\pi^*\colon E(\mathbf{BG}_m) \to E(\mathfrak{L})$ is an equivalence by Definition~\ref{def:axiomatic-qsm-class}.)  
\end{definition}

\begin{definition}
\label{def:universal-mgl}
We define the algebraic cobordism $\MGL_\mathcal{C}$ as the \emph{universal oriented spectrum} on $\mathcal{C}$. Thus, $\MGL_\mathcal{C}$ is an $\mathbb{E}_\infty$-ring spectrum that is equipped with a universal orientation class
\[
    c_1^{\mathrm{univ}} \in \MGL_\mathcal{C}(\mathbf{BG}_m),
\]
such that for any oriented $\mathbb{E}_\infty$-ring spectrum $(E, c_1^E)$ on $\mathcal{C}$, there exists a unique morphism of $\mathbb{E}_\infty$-ring spectra
\[
    \varphi_E\colon \MGL_\mathcal{C} \longrightarrow E
\]
that satisfies $\varphi_E(c_1^{\mathrm{univ}}) \simeq c_1^E$. 
\end{definition}

\begin{definition}
\label{scheme-setting}
Suppose that $\mathcal{C}$ is an $\infty$-category arising from schemes (or
derived schemes) over a base scheme $S$. For a morphism
$f\colon X\to Y$ in $\mathcal{C}$, we say that  $f$ is {\it quasi-smooth} if
its cotangent complex $\mathbb{L}_{X/Y}$ is perfect of Tor-amplitude in
$[-1,0]$. In the classical scheme context, this is equivalent to the locally complete intersection (l.c.i.) condition. Let
\[
\mathrm{QSm}(\mathcal{C})\subset\mathrm{Mor}(\mathcal{C}).
\]
denote the class of quasi-smooth morphisms. 
\end{definition}
\begin{proposition}
\label{prop:qsm-admissible}
In the case of Definition~\ref{scheme-setting}, the class
$\mathbf{S}=\mathrm{QSm}(\mathcal{C})$ satisfies the conditions (1), (2), and (3) outlined in Definition~\ref{def:axiomatic-qsm-class}.
\end{proposition}

\begin{proof}
For Condition~(1), quasi-smooth morphisms are stable under equivalences, composition through the transitivity sequence of the cotangent complex, and arbitrary derived base change by~\cite[p.~1270, Proposition 7.2.4.23 and p.~1300, Proposition 7.3.3.7]{HA}.

For Condition~(2), the codiagonal $\nabla\colon X\amalg X\to X$ restricts on each summand to the identity $\mathrm{id}_X\colon X\to X$, which is an isomorphism. In particular, $\nabla$ is locally \'etale. Therefore, it is quasi-smooth. 

For Condition~(3), the line bundle projection $\pi\colon \mathfrak{L}\to \mathbf{BG}_m$ is smooth, and its zero-section is a locally complete intersection; hence, both morphisms are quasi-smooth. Since $\pi$ is a line-bundle projection, it is an $\mathbb{A}^1$-weak equivalence (see \cite{Hoyois2017SixOps} or \cite{Jardine2}). \qed
\end{proof}

By Proposition~\ref{prop:qsm-admissible}, the class $\mathrm{QSm}(\mathcal{C})$ is an admissible class, so the span-theoretic algebraic cobordism $\MGL_{\mathrm{H}(S)}$ is well-defined. By Definition~\ref{def:universal-mgl}, $\MGL_{\mathrm{H}(S)}$ is initial in the $\infty$-category of oriented $\mathbb{E}_\infty$-ring spectra on $\mathrm{Span}_{\mathrm{QSm}(\mathrm{H}(S))}(\mathrm{H}(S))$.
\begin{proposition}
\label{prop:gesn-comparison}
Let $S$ be a scheme, and let $\mathrm{H}(S)$ be the $\infty$-category of motivic spaces introduced in~\cite{MV}.
Then the algebraic cobordism $\MGL_{\mathrm{H}(S)}$ produced by the span formalism with admissible class given by quasi-smooth morphisms is canonically equivalent to Voevodsky's algebraic cobordism spectrum $\MGL$.
\end{proposition}
\begin{proof}
By \cite[Theorem~2.7]{zbMATH05367298}, 
Voevodsky's algebraic cobordism spectrum $\MGL$~\cite{V}, identified with the motivic Thom-spectrum by Gepner--Snaith~\cite[Definition~3.1]{GeSn}, serves as an initial object of oriented motivic $\mathbb{E}_\infty$-ring spectra.

In the context of motivic spectra, our span formalism is generated exclusively by usual pull-backs and oriented quasi-smooth Gysin operations. Thus, the universal properties of $\MGL_{\mathrm{H}(S)}$ as presented in Definition~\ref{def:universal-mgl} and of $\MGL$, which belong to the same $\infty$-category of oriented motivic theories, imply their weak equivalence. \qed
\end{proof}

\section{The projective bundle formula and Chern class formalism}
\label{sec:pbf}
In this section, we establish the projective bundle formula and develop the Chern-class formalism, including the Whitney sum formula. Throughout, we work under the following standing assumptions about the ambient underlying $\infty$-category $\mathcal{C}$.

\begin{assumption}\label{assum:standing_section4}
We assume that $\mathcal{C}$ is locally presentable (so that the required quotient objects exist) and satisfies the following conditions:
\begin{enumerate}
    \item \textbf{Group structure:} The tensor-product monoid structure on $\mathbf{BG}_m$ is group-like, making the associated $\mathbb{G}_m$ a group object.
    \item \textbf{Fiber identification:} The fiber of $\mathfrak{L}$ over a given rational point $* \to \mathbf{BG}_m$ is identified with $\mathbb{A}^1$.
    \item \textbf{Equivalence and Torsor property:} Base change along the projection $\mathfrak{L}^\circ \to \mathbf{BG}_m$ induces a canonical equivalence 
    \[
    \Phi \colon \mathfrak{L}^\circ \times \mathbb{A}^1 \to \mathfrak{L}^\circ \times_{\mathbf{BG}_m} \mathfrak{L}.
    \]
    Furthermore, the restriction $\Phi^\circ$ makes $\mathfrak{L}^\circ\to\mathbf{BG}_m$ a principal $\mathbb{G}_m$-torsor in the sense of Definition~\ref{def:tautological_shear}.
    \item \textbf{Base-change stability:} The base changes of this principal $\mathbb{G}_m$-torsor used in our subsequent constructions remain principal $\mathbb{G}_m$-torsors.
\end{enumerate}
\end{assumption}

\begin{remark}
Condition (4) in Assumption~\ref{assum:standing_section4} requires some categorical explanation. It is automatic if $\mathcal{C}$ is an $\infty$-topos, because effective epimorphisms are strictly stable under base change. However, in general, after the $\mathbb{A}^1$-localization of an $\infty$-category, the localization functor fails to be left exact. 

Nevertheless, this technical problem disappears once we pass to a spectrum-valued sheaf $E$. The $\mathbb{A}^1$-localization of spectrum-valued sheaves constitutes a Bousfield localization of a \emph{stable} $\infty$-category, indicating that its colimit-preserving localization functor is exact and preserves finite limits. Thus, the necessary torsor descent can be verified after applying $E$ in this stable sheaf category. To maintain geometric flexibility, we incorporate torsor base-change stability into the standing assumptions.
\end{remark}

\subsection{The $\infty$-categorical shear action}
This subsection constructs the $\mathbb{G}_m$-action on the
punctured universal bundle utilized in the cellular construction.
Let
\[
\pi^\circ=\pi\circ i\colon \mathfrak{L}^\circ\to \mathbf{BG}_m
\]
denote the restriction of the universal line bundle to the punctured line
object.

\begin{definition}[The $\infty$-categorical $\mathbb{G}_m$-action]
\label{def:tautological_shear}
The pull-back of the universal line bundle $\pi\colon \mathfrak{L}\to \mathbf{BG}_m$
along $\pi^\circ \colon \mathfrak{L}^\circ \to \mathbf{BG}_m$ induces a canonical equivalence
\[
\Phi \colon \mathfrak{L}^\circ \times \mathbb{A}^1 \to \mathfrak{L}^\circ \times_{\mathbf{BG}_m} \mathfrak{L}
\]
over $\mathfrak{L}^\circ$. Restricting $\Phi$ to the open subspace
$\mathbb{G}_m \subset \mathbb{A}^1$ removes the zero-section and induces an
equivalence
\[
\Phi^\circ\colon \mathfrak{L}^\circ\times \mathbb{G}_m \to \mathfrak{L}^\circ\times_{\mathbf{BG}_m}\mathfrak{L}^\circ.
\]
We define the canonical map $\mu \colon \mathfrak{L}^\circ \times \mathbb{G}_m \to \mathfrak{L}^\circ$ as the composition $\mu = \mathrm{pr}_2 \circ \Phi^\circ$, where $\mathrm{pr}_2$ is the second projection of the fiber product. By the construction, this formal trivialization $\Phi^\circ$ ensures that $\mu$ defines a principal $\mathbb{G}_m$-action on $\mathfrak{L}^\circ$ without choosing a global rational base point.

Here and below, a principal $\mathbb{G}_m$-torsor over an object $X$ means a
morphism $P \to X$ equipped with a $\mathbb{G}_m$-action over $X$ such that
the shear map
\[
P\times\mathbb{G}_m\longrightarrow P\times_X P
\]
is an equivalence and $P\to X$ is an effective epimorphism, so that $X$ is
recovered as the quotient of the action groupoid. For
$\pi^\circ\colon\mathfrak{L}^\circ\to\mathbf{BG}_m$, the shear map is exactly
$\Phi^\circ$.
\end{definition}

\begin{remark}[The classical motivic model]
\label{rem:contractibility-Lcirc}
In the usual motivic category of schemes, it is classical that
$\mathbb{A}^\infty\setminus\{0\}$ is $\mathbb{A}^1$-contractible. Indeed,
we compare the standard inclusions $x\mapsto (x,0)$ with the shifted inclusions
$x\mapsto (x,1)$. These two systems are connected by the evident
$\mathbb{A}^1$-homotopy $(x,t)$, while the shifted punctured system interleaves
with the filtered system of affine spaces $\mathbb{A}^n$, whose
colimit is $\mathbb{A}^1$-contractible.
\end{remark}

In $\mathbb{A}^1$-homotopy theory, there is a classical identification
$\mathbf{BG}_m\simeq \mathbb{P}^\infty$ in the $\mathbb{A}^1$-homotopy category;
equivalently, the universal torsor $\mathfrak{L}^\circ\to \mathbf{BG}_m$ is modeled by
$\mathbb{A}^\infty\setminus\{0\}\to \mathbb{P}^\infty$.
In this model, the total space $\mathfrak{L}^\circ \simeq \mathbb{A}^\infty\setminus\{0\}$ is
$\mathbb{A}^1$-contractible, justifying its role as the universal space
$E\mathbb{G}_m$. By Remark~\ref{rem:contractibility-Lcirc}, this classical construction
can be presented using the two rational points $0,1\colon *\to\mathbb{A}^1$ to describe
the bonding of the $\mathbb{A}^1$-cells. 

In this model, when $\mathbb{G}_m$ is a group object and $\mathfrak{L}^\circ$ is contractible,
the principal fibration $\mathbb{G}_m \to \mathfrak{L}^\circ \to \mathbf{BG}_m$
induces an equivalence of group objects
$\mathbb{G}_m \simeq \Omega \mathbf{BG}_m \simeq \Omega B \Gm$.
Delooping provides the canonical identification
$B\mathbb{G}_m \simeq \mathbf{BG}_m$, ensuring that the quotient of
punctured total spaces under the tautological action produces the expected
projective bundles.

In our setting, the formal arguments presented below do not directly use this
contractibility. The action $\mu$ is
solely induced by the equivalence $\Phi$ from Definition~\ref{def:tautological_shear} and does not require any global choice
of a base point. The restriction of $\Phi$ identifies the punctured
pull-back of the universal line bundle with the trivial 
$\mathbb{G}_m$-torsor over $\mathfrak{L}^\circ$.

\subsection{Definition of the cellular pushout square}

This subsection defines the cube used in the cellular construction and
records the resulting push-out square.
Whenever an object from the $(n-1)$-stage is viewed over $BT_n$, we
implicitly pull it back along the projection
$q_n\colon BT_n\to BT_{n-1}$ while forgetting the last factor, denoting the
pull-back by the same symbol.

\begin{definition}\label{def:pbf-cube}
For each $n\ge 1$, write $BT_n=(\mathbf{BG}_m)^n$. For each $1 \le i \le n$, let $\mathrm{pr}_i\colon BT_n\to \mathbf{BG}_m$ denote the $i$-th projection, and set 
\[
\mathfrak{L}_i=\operatorname{pr}_i^*(\mathfrak{L})\to BT_n, \quad \mathfrak{L}_i^\circ = \mathrm{pr}_i^*(\mathfrak{L}^\circ), 
\]
and
\[
\mathfrak{L}^n=\underbrace{\mathfrak{L}_1\times_{BT_n}\cdots\times_{BT_n}\mathfrak{L}_n}_{n\ \text{factors}},
\qquad
\mathcal U_i=\mathfrak{L}_1\times_{BT_n}\cdots\times_{BT_n}\mathfrak{L}_i^\circ\times_{BT_n}\cdots\times_{BT_n}\mathfrak{L}_n
\ \ (1\le i\le n),
\]
where in $\mathcal U_i$ the $i$-th factor is $\mathfrak{L}_i^\circ$ and the other factors are $\mathfrak{L}_j$.
The punctured total space is defined by
\[
(\mathfrak{L}^n)^\circ=\bigcup_{i=1}^n \mathcal U_i\subset \mathfrak{L}^n. 
\]
We define the projective bundle $\mathbb{P}(\mathfrak{L}^n)$ as the quotient of $(\mathfrak{L}^n)^\circ$ by the $\mathbb{G}_m$-action of Definition~\ref{def:tautological_shear}.

For any nonempty subset $I\subseteq\{1,\dots,n\}$, define
\[
\mathcal U_I=\lim_{i\in I}(\mathcal U_i\to \mathfrak{L}^n)
\simeq \prod_{i\in I,\,\mathfrak{L}^n}\mathcal U_i\subset (\mathfrak{L}^n)^\circ.
\]
The assignment $I\mapsto \mathcal U_I$ forms an $n$-cube by inclusions. The quotients $\mathcal U_i/\mathbb{G}_m$ form the standard open cover of $\mathbb{P}(\mathfrak{L}^n)$, so $\mathbb{P}(\mathfrak{L}^n)$ is identified with the colimit of the punctured quotient diagram:
\[
\mathbb{P}(\mathfrak{L}^n)\simeq \operatorname*{colim}_{\emptyset\neq I\subseteq\{1,\dots,n\}}\bigl(\mathcal U_I/\mathbb{G}_m\bigr).
\]
For $n\ge 2$ and $1\le i\le n-1$, consider the quotient
\[
\mathfrak{W}_{i,n}=\bigl(\mathfrak{L}_i\times_{BT_n}\mathfrak{L}_n^\circ\bigr)/\mathbb{G}_m
\longrightarrow \mathfrak{L}_n^\circ/\mathbb{G}_m,
\]
where $\mathbb{G}_m$ acts diagonally via the tautological action on
$\mathfrak{L}_n^\circ$ and fiberwise scalar multiplication on
$\mathfrak{L}_i$. The morphism $\mathfrak{L}_n^\circ\to BT_n$ is the
base change of $\pi^\circ\colon\mathfrak{L}^\circ\to\mathbf{BG}_m$ along
$\mathrm{pr}_n\colon BT_n\to\mathbf{BG}_m$; hence, by the standing
assumptions of this section, it is a principal $\mathbb{G}_m$-torsor. This
equivariant line bundle over $\mathfrak{L}_n^\circ$ therefore descends to a
line bundle on $BT_n$. Let
\[
\chi_{i,n}\colon BT_n=(\mathbf{BG}_m)^n\longrightarrow \mathbf{BG}_m
\]
denote the classifying map of $\mathfrak{W}_{i,n}\to BT_n$. Then
\[
\mathfrak{W}_{i,n}\simeq
BT_n\times_{\chi_{i,n},\mathbf{BG}_m}\mathfrak{L}.
\]
Set
\[
\mathfrak{W}_n=\underbrace{\mathfrak{W}_{1,n}\times_{BT_n}\cdots\times_{BT_n}\mathfrak{W}_{n-1,n}}_{n-1\ \text{factors}},
\]
and let $\mathfrak{W}_n^\circ\subset \mathfrak{W}_n$ denote its punctured
total space. The quotient descriptions identify the $n$-th affine chart as
\[
\mathfrak{W}_n\simeq \mathcal U_n/\mathbb{G}_m.
\]
For each nonempty subset $J\subseteq\{1,\dots,n-1\}$, the intersection of
the standard charts of $\mathfrak{W}_n$ indexed by $J$ identifies
with $\mathcal U_{J\cup\{n\}}/\mathbb{G}_m$. Hence, descent for the standard
open cover of the punctured total space gives
\[
\mathfrak{W}_n^\circ\simeq
\operatorname*{colim}_{\emptyset\neq J\subseteq\{1,\dots,n-1\}}
\bigl(\mathcal U_{J\cup\{n\}}/\mathbb{G}_m\bigr).
\]
Geometrically, $\mathfrak{W}_n^\circ$ is the locus where the $n$-th
coordinate is non-zero and at least one of the first $(n-1)$-coordinates is
non-zero.
\end{definition}

By separating the $n$-th chart from the union of the first $(n-1)$ charts in the cubical cover, we obtain the following compactification of the affine chart $\mathfrak{W}_n$ using the hyperplane $\mathbb{P}(\mathfrak{L}^{n-1})$.

\begin{lemma}[The cellular pushout square]\label{lem:cell-square-derived}
For each $n \geq 2$, there exists a canonical homotopy coCartesian square in
$\mathcal{C}$:
\[
\xymatrix@C=4pc@R=3pc{
    \mathfrak{W}_n^\circ  \ar[r]^-{i_n^{\Box}}  \ar[d]_{p_n} & \mathfrak{W}_n \ar[d]^{j_n} \\
   \mathbb{P}(\mathfrak{L}^{n-1})  \ar[r]_-{i_n} & \mathbb{P}(\mathfrak{L}^n)
}
\]
Here, $i_n^{\Box}$ denotes the punctured-total-space inclusion,
$p_n\colon \mathfrak{W}_n^\circ\to \mathbb{P}(\mathfrak{L}^{n-1})$ represents the
quotient map, $j_n\colon \mathfrak{W}_n\to \mathbb{P}(\mathfrak{L}^n)$ corresponds
to the open immersion of the distinguished chart associated with the
$n$-th summand, and $i_n$ indicates the hyperplane inclusion induced by the first
$(n-1)$ summands.
\end{lemma}

\begin{proof}
By Definition~\ref{def:pbf-cube}, one has
\[
\mathbb{P}(\mathfrak{L}^n)\simeq
\operatorname*{colim}_{\emptyset\neq I\subseteq\{1,\dots,n\}}
\bigl(\mathcal U_I/\mathbb{G}_m\bigr).
\]
Restricting to the subsets $I\subseteq\{1,\dots,n-1\}$ recovers the same
punctured quotient diagram for $\mathbb{P}(\mathfrak{L}^{n-1})$; hence, 
\[
\mathbb{P}(\mathfrak{L}^{n-1})\simeq
\operatorname*{colim}_{\emptyset\neq I\subseteq\{1,\dots,n-1\}}
\bigl(\mathcal U_I/\mathbb{G}_m\bigr).
\]
The remaining chart is
$\mathfrak{W}_n\simeq \mathcal U_n/\mathbb{G}_m$, and its overlap with the
preceding colimit is the colimit of the mixed intersections
$\mathcal U_{J\cup\{n\}}/\mathbb{G}_m$ for
$\emptyset\neq J\subseteq\{1,\dots,n-1\}$. By
Definition~\ref{def:pbf-cube}, this intersection is the punctured
total space $\mathfrak{W}_n^\circ$. Therefore, the projective bundle $\mathbb{P}(\mathfrak{L}^n)$
is obtained by gluing $\mathfrak{W}_n$ to
$\mathbb{P}(\mathfrak{L}^{n-1})$ along $\mathfrak{W}_n^\circ$. The square mentioned in the assertion is homotopy coCartesian. \qed
\end{proof}

\begin{lemma}
\label{lem:split-thom}
For each $n\ge 2$, the bundle $\mathfrak{W}_n\to BT_n$ is a split
rank-$(n-1)$ vector bundle with line-bundle factors
$\mathfrak{W}_{i,n}\to BT_n$ $(1\le i\le n-1)$. There exists a canonical equivalence
\[
\mathrm{Th}(\mathfrak{W}_n)\simeq
\mathrm{Th}(\mathfrak{W}_{1,n})\wedge\cdots\wedge\mathrm{Th}(\mathfrak{W}_{n-1,n}).
\]
If $E$ is an oriented cohomology theory and $\tau_n\in E(\mathrm{Th}(\mathfrak{W}_n))$
denotes the external product of the Thom classes of the factors
$\mathfrak{W}_{i,n}$, then multiplication by $\tau_n$ induces an
equivalence of $E(BT_n)$-module spectra
\[
E(BT_n)\to E(\mathrm{Th}(\mathfrak{W}_n)).
\]
In particular, $E(\mathrm{Th}(\mathfrak{W}_n))$ is a free rank-one
$E(BT_n)$-module spectrum.
\end{lemma}

\begin{proof}
By the standing assumptions of this section,
$\mathfrak{L}_n^\circ \to BT_n$ is a principal $\mathbb{G}_m$-torsor. For
each $1\le i\le n-1$, the induced map $\mathfrak{W}_{i,n}\to BT_n$ is a line bundle. Taking their
fiber product over $BT_n$ gives the split rank-$(n-1)$ vector bundle
$\mathfrak{W}_n\to BT_n$.

The Thom space of a split vector bundle is obtained from the iterated
pushout product of the punctured zero-section inclusions of its line-bundle
factors. Hence, the first equivalence follows by iterating
Lemma~\ref{lem:thom-symmetric-monoidal}. For each factor, since
$\mathfrak{W}_{i,n}$ is the pull-back of the universal line bundle along
$\chi_{i,n}\colon BT_n\to\mathbf{BG}_m$, the pull-back of the orientation equivalence
\[
E(\mathbf{BG}_m)\to E(\mathrm{Th}(\mathfrak{L}))
\]
along $\chi_{i,n}$ induces an equivalence of
$E(BT_n)$-module spectra
\[
E(BT_n)\to E(\mathrm{Th}(\mathfrak{W}_{i,n})).
\]
After taking the tensor product of these $(n-1)$ equivalences, $E(\mathrm{Th}(\mathfrak{W}_n))$ can be identified with a free rank-one $E(BT_n)$-module generated by the Thom class $\tau_n$ of the product.
\qed
\end{proof}

\subsection{The projective bundle formula}

We establishes the projective bundle formula for the split
universal bundle at the level of $E$-module spectra. Fix an oriented
cohomology theory $E$ on $\mathcal{C}$. 
\begin{theorem}
\label{thm:pbf-universal}
Under the standing assumptions of this section, the canonical morphism
\[
\Phi_n \colon \bigoplus_{i=0}^{n-1} E(BT_n)\cdot \xi^i \longrightarrow E(\mathbb{P}(\mathfrak{L}^n))
\]
is a weak equivalence of $E(BT_n)$-module spectra, where
$\xi=\xi_{\mathfrak{L}^n}$ is the hyperplane class on
$\mathbb{P}(\mathfrak{L}^n)$. Therefore, after applying $\pi_{-*}$,
it induces an isomorphism
\[
\bigoplus_{i=0}^{n-1} E^*(BT_n)\cdot \xi^i \to E^*(\mathbb{P}(\mathfrak{L}^n))
\]
of $E^*(BT_n)$-modules.
\end{theorem}

\begin{remark}
\label{rem:top-thom-normalization}
Since $
i_n^*(\xi_{\mathfrak{L}^n})=\xi_{\mathfrak{L}^{n-1}}$,
the lower-power summand $
\bigoplus_{i=0}^{n-2}E(BT_n)\cdot\xi_{\mathfrak{L}^n}^i$ is identified via $i_n^*$ with $E(\mathbb{P}(\mathfrak{L}^{n-1}))$ by the
inductive hypothesis. Thus, the fiber term $E(\mathrm{Th}(\mathfrak{W}_n))$
represents the residual rank-one quotient of $E(\mathbb{P}(\mathfrak{L}^n))$
modulo the lower powers; thus, stating that its generator corresponds to
$\xi_{\mathfrak{L}^n}^{n-1}$ is understood within that quotient.
\end{remark}

\begin{proof}
We prove this by induction on $n$.
For $n=1$, the assertion is obvious. For $n\ge 2$, let
\[
q_n\colon BT_n=(\mathbf{BG}_m)^n\longrightarrow BT_{n-1}=(\mathbf{BG}_m)^{n-1}
\]
be the projection forgetting the last factor. By the convention fixed at the
start of this subsection, $\mathbb{P}(\mathfrak{L}^{n-1})$ denotes its
pull-back along $q_n$, and its hyperplane class pulls back to
$\xi_{\mathfrak{L}^n}$.

By Definition~\ref{def:oriented-cohomology}, $E$ is excisive, and the cellular push-out square of Lemma~\ref{lem:cell-square-derived}
gives a fiber sequence of $E(BT_n)$-module spectra
\[
E(\mathbb{P}(\mathfrak{L}^{n-1}))
\xleftarrow{i_n^*}
E(\mathbb{P}(\mathfrak{L}^{n}))
\leftarrow
E(\mathrm{Th}(\mathfrak{W}_n)).
\]

By the assumption of induction, the composite
\[
\bigoplus_{i=0}^{n-2}E(BT_n)\cdot\xi_{\mathfrak{L}^n}^i
\longrightarrow E(\mathbb{P}(\mathfrak{L}^n))
\xrightarrow{i_n^*}
E(\mathbb{P}(\mathfrak{L}^{n-1}))
\]
is an equivalence of $E(BT_n)$-module spectra. Therefore, the lower powers
provide a section of $i_n^*$, and the fiber sequence splits as
\[
E(\mathbb{P}(\mathfrak{L}^n))
\simeq
E(\mathbb{P}(\mathfrak{L}^{n-1})) \oplus E(\mathrm{Th}(\mathfrak{W}_n)).
\]
By Lemma~\ref{lem:split-thom}, the second summand is a free rank-one
$E(BT_n)$-module spectrum. By Remark~\ref{rem:top-thom-normalization}, its
generator is represented modulo the lower-power summand by the top power
$\xi_{\mathfrak{L}^n}^{n-1}$. Replacing the chosen generator by the actual
class $\xi_{\mathfrak{L}^n}^{n-1}$ corresponds to an upper-triangular change of
basis with an identity diagonal, thus providing an equivalence
\[
E(\mathbb{P}(\mathfrak{L}^n))
\simeq
E(\mathbb{P}(\mathfrak{L}^{n-1})) \oplus E(BT_n)\cdot \xi_{\mathfrak{L}^n}^{n-1}.
\]

The assumption of induction implies that the case of $n-1$ along $q_n$ gives an
equivalence of $E(BT_n)$-module spectra
\[
E(\mathbb{P}(\mathfrak{L}^{n-1}))
\simeq
\bigoplus_{i=0}^{n-2}E(BT_n)\cdot\xi_{\mathfrak{L}^n}^i.
\]
Thus, we obtain the equivalence
\[
E(\mathbb{P}(\mathfrak{L}^n))
\simeq
\bigoplus_{i=0}^{n-1}E(BT_n)\cdot\xi_{\mathfrak{L}^n}^i,
\]
which is the claimed weak equivalence. The final statement on homotopy groups
follows by applying $\pi_{-*}$. \qed
\end{proof}

\subsection{From the split universal bundle to arbitrary vector bundles}

We first transition from the split rank-$n$ bundle
$\mathfrak{L}^n$ over $BT_n=(\mathbf{BG}_m)^n$ to the universal rank-$n$
bundle $\mathfrak{V}_n$ over $\mathbf{BGL}_n$ via the $\Sigma_n$-comparison. Then,
we deduce the result for an arbitrary rank-$n$ vector bundle using the
pull-back along its classifying map. Only at the end do we apply
$\pi_{-*}$ to recover the usual graded projective bundle formula.

\begin{definition}[The Grassmann objects and the splitting principle property]
\label{prop:bgl-assumption}
For each $n\ge 1$, let $\mathbf{BGL}_n$ denote the {\it Grassmann object of rank $n$},
and let
\[
\pi_n\colon \mathfrak{V}_n \to \mathbf{BGL}_n
\]
be the universal rank-$n$ vector bundle whose projection belongs to $\mathbf{S}$.
Moreover, any rank-$n$ vector bundle $V\to X$ is classified by a map
$f_V\colon X\to \mathbf{BGL}_n$ with
$V\simeq f_V^*(\mathfrak{V}_n)$.
We say that $E$ has the {\it splitting principle property} if, for the map $\mathrm{u}_n\colon BT_n\to \mathbf{BGL}_n$, the
pull-back along $\mathrm{u}_n$ induces weak equivalences
\[
E(\mathbf{BGL}_n)\to E(BT_n)^{\Sigma_n},
\qquad
E(\mathbb{P}(\mathfrak{V}_n)) \to
E(\mathbb{P}(\mathfrak{L}^n))^{\Sigma_n}.
\]
Here, the second map uses the canonical identification
$\mathbb{P}(\mathrm{u}_n^*\mathfrak{V}_n)\simeq\mathbb{P}(\mathfrak{L}^n)$.
These equivalences are compatible with the module structures over the first
equivalence and carry the hyperplane class $\xi_{\mathfrak{V}_n}$ to the
fixed hyperplane class $\xi_{\mathfrak{L}^n}$.
We require the same comparison after finite products of Grassmann objects.
In particular, the pull-back along
$\mathrm{u}_{n_1}\times\cdots\times\mathrm{u}_{n_r}$ identifies
$E(\prod_a\mathbf{BGL}_{n_a})$ with
$E(\prod_aBT_{n_a})^{\prod_a\Sigma_{n_a}}$ compatibly with the universal
bundles on each factor. Finally, this universal comparison is assumed to be
compatible with classifying pull-backs: for every map
$f\colon X\to\mathbf{BGL}_n$ classifying a rank-$n$ vector bundle
$V=f^*\mathfrak{V}_n$, the canonical extension-of-scalars map
\[
E(X)\otimes_{E(\mathbf{BGL}_n)}E(\mathbb{P}(\mathfrak{V}_n))
\longrightarrow E(\mathbb{P}(V))
\]
is a weak equivalence that carries the universal hyperplane class to
$\xi_V$.
\end{definition}

\begin{proposition}
\label{prop:pbf-bgl}
Under the assumptions of this section, we further assume that 
$\mathcal{C}$ contains Grassmann objects of all ranks. For any oriented
cohomology theory $E$ that exhibits the splitting principle property, the canonical morphism
\[
\Phi_{\mathfrak{V}_n}\colon \bigoplus_{i=0}^{n-1}E(\mathbf{BGL}_n)\cdot \xi^i
\longrightarrow E(\mathbb{P}(\mathfrak{V}_n))
\]
is a weak equivalence of $E(\mathbf{BGL}_n)$-module spectra, where
$\xi=\xi_{\mathfrak{V}_n}$ represents the hyperplane class on
$\mathbb{P}(\mathfrak{V}_n)$.
\end{proposition}

\begin{proof}
By Theorem~\ref{thm:pbf-universal}, the split rank-$n$ bundle
$\mathfrak{L}^n$ over $BT_n$ satisfies a weak equivalence
\[
\Phi_n\colon \bigoplus_{i=0}^{n-1}E(BT_n)\cdot \xi_{\mathfrak{L}^n}^i
\to
E(\mathbb{P}(\mathfrak{L}^n)).
\]
The permutation action of $\Sigma_n$ on the universal line factors of
$\mathfrak{L}^n$ induces $\Sigma_n$-actions on both sides, making
$\Phi_n$ is equivariant. The hyperplane class $\xi_{\mathfrak{L}^n}$ is
$\Sigma_n$-invariant. Since fixed points in spectra are limits and
$\Sp$ is stable, taking $\Sigma_n$-fixed points commutes with this finite
direct sum. Thus, the passing to fixed points induces a weak equivalence of
$E(BT_n)^{\Sigma_n}$-module spectra
\[
\bigoplus_{i=0}^{n-1}E(BT_n)^{\Sigma_n}\cdot \xi_{\mathfrak{L}^n}^i
\to
E(\mathbb{P}(\mathfrak{L}^n))^{\Sigma_n}.
\]
According to the comparison equivalences in Definition~\ref{prop:bgl-assumption},
these fixed-point spectra correspond to $E(\mathbf{BGL}_n)$ and
$E(\mathbb{P}(\mathfrak{V}_n))$ in a way compatible with the module
structures and the hyperplane classes
$\xi_{\mathfrak{V}_n}\mapsto\xi_{\mathfrak{L}^n}$. Transporting the
preceding equivalence along these identifications induces a weak equivalence
\[
\Phi_{\mathfrak{V}_n}\colon \bigoplus_{i=0}^{n-1}E(\mathbf{BGL}_n)\cdot \xi^i
\to
E(\mathbb{P}(\mathfrak{V}_n)).
\] \qed
\end{proof}

\begin{corollary}
\label{cor:pbf-general}
Under the standing assumptions of this section, we further assume that 
$\mathcal{C}$ contains Grassmann objects of all ranks. For any oriented
cohomology theory $E$ that possesses the splitting principle property and for any rank-$n$ vector bundle $V \to X$ over $X \in \mathcal{C}$, the
canonical morphism
\[
    \Phi_V \colon \bigoplus_{i=0}^{n-1} E(X)\cdot \xi^i \longrightarrow E(\mathbb{P}(V))
\]
is a weak equivalence of $E(X)$-module spectra, where $\xi=\xi_V$ is the
hyperplane class on $\mathbb{P}(V)$. After applying
$\pi_{-*}$, it induces an isomorphism
\[
\bigoplus_{i=0}^{n-1}E^*(X)\cdot\xi^i \xrightarrow{\ \cong\ } E^*(\mathbb{P}(V))
\]
of $E^*(X)$-modules.
\end{corollary}

\begin{proof}
Let $f_V\colon X\to \mathbf{BGL}_n$ denote the classifying map of $V \to X$. The projection $V \to X$ commutes with
pull-backs and provides a canonical equivalence
\[
\mathbb{P}(V)\simeq \mathbb{P}(f_V^*\mathfrak{V}_n)
\simeq X\times_{\mathbf{BGL}_n}\mathbb{P}(\mathfrak{V}_n).
\]
Under this identification, the hyperplane class $\xi_V$ is the pull-back of
$\xi_{\mathfrak{V}_n}$. Therefore, the compatibility of extension-of-scalars in 
Definition~\ref{prop:bgl-assumption} identifies the pull-back of
$\Phi_{\mathfrak{V}_n}$ along $f_V$ with
\[
\Phi_V\colon \bigoplus_{i=0}^{n-1}E(X)\cdot \xi^i
\longrightarrow E(\mathbb{P}(V)).
\]
Therefore, Proposition~\ref{prop:pbf-bgl} implies that $\Phi_V$ is a weak
equivalence of $E(X)$-module spectra. The final statement follows from
applying $\pi_{-*}$. \qed
\end{proof}

\begin{remark}
Write $E^*=\pi_{-*}E$. The projective bundle formula implies that
\[
E^*(BT_n) \cong E^*(pt)[[x_1, \dots, x_n]],
\]
where $BT_n \simeq (\mathbf{BG}_m)^n$ and $x_1,\dots,x_n$ are the universal
Chern roots. The symmetric group $\Sigma_n$ acts on $BT_n$ by permuting the
universal line bundles, and therefore, it acts on $E^*(BT_n)$ by permuting the roots. Under
this action, the universal Chern classes $c_k$ correspond to the elementary symmetric
polynomials in $x_1,\dots,x_n$, which leads to the standard
identification
$E^*(BT_n)^{\Sigma_n}\simeq E^*(\mathbf{BGL}_n)$.
\end{remark}

\subsection{Chern classes and Whitney sums}

This subsection presents the consequences of the projective
bundle formula for a fixed oriented theory $E$ that satisfies the splitting principle property. 

\begin{definition}
Let $V \to X$ be a rank-$n$ vector bundle, and let $\pi\colon \mathbb{P}(V) \to X$ be the projective bundle. By the projective bundle formula (Corollary~\ref{cor:pbf-general}), $E(\mathbb{P}(V))$ is a free module over the $\mathbb{E}_\infty$-ring $E(X)$ with basis $1, \xi, \dots, \xi^{n-1}$, where $\xi$ represents the first Chern class of the universal line bundle. Since $\xi^n \in E(\mathbb{P}(V))$ must be uniquely expressed as an $E(X)$-linear combination of this basis, there exist unique classes $c_k(V) \in \pi_0 E(X)$ for $k=1, \dots, n$, called the \emph{Chern classes} of $V$, such that the relation
\[
    \sum_{i=0}^n (-1)^i \pi^*(c_i(V)) \xi^{n-i} = 0
\]
holds in $E(\mathbb{P}(V))$, where we set $c_0(V) = 1$.
\end{definition}

\begin{theorem}
\label{thm:whitney-sum}
Under the standing assumptions of this section, we further assume that 
$\mathcal{C}$ contains Grassmann objects of all ranks. For any oriented
cohomology theory $E$ exhibiting the splitting principle property, and for two vector bundles $V$ and $W$ over $X$, the total Chern class $c(V \oplus W) = \sum_{k} c_k(V \oplus W)$ satisfies
\[
    c(V \oplus W) = c(V) c(W).
\]
In terms of components, $c_k(V \oplus W) = \sum_{i+j=k} c_i(V) c_j(W)$.
\end{theorem}

\begin{proof}
This is the usual splitting-principle proof of the Whitney sum formula (see Fulton~\cite[Remark~3.2.3]{FultonIT}). We only explain how the classical argument is implemented under
Definition~\ref{prop:bgl-assumption}.

Let $m$ be the rank of $V$, $n$ be the rank of $W$, and let
$p_1:\mathbf{BGL}_m\times \mathbf{BGL}_n\to \mathbf{BGL}_m$ and $p_2: \mathbf{BGL}_m\times \mathbf{BGL}_n \to \mathbf{BGL}_n$ denote the projections. Since the bundle
$p_1^*\mathfrak{V}_m\oplus p_2^*\mathfrak{V}_n$ has rank $m+n$, the
universal property of the Grassmann object of rank $m+n$ provides a classifying map
$\mu_{m,n}\colon \mathbf{BGL}_m\times \mathbf{BGL}_n\to \mathbf{BGL}_{m+n}$
classifying it. Since the Chern classes defined above are natural under
pull-back, it is enough to prove the universal identity
\[
    \mu_{m,n}^*c(\mathfrak{V}_{m+n})
    = p_1^*c(\mathfrak{V}_m)\,p_2^*c(\mathfrak{V}_n)
\]
in $\pi_0 E(\mathbf{BGL}_m\times \mathbf{BGL}_n)$, since any pair
$(V,W)$ is obtained from $(\mathfrak{V}_m,\mathfrak{V}_n)$ by pull-back
along the joint classifying map
$g=(f_V,f_W)\colon X\to \mathbf{BGL}_m\times \mathbf{BGL}_n$.

After pull-back along the maximal-torus map
$\mathrm{u}_m\times \mathrm{u}_n\colon BT_m\times BT_n\to
\mathbf{BGL}_m\times \mathbf{BGL}_n$, the relevant universal bundle is split as a direct sum of
universal line bundles. Writing $x_1,\dots,x_m$ and $y_1,\dots,y_n$ for the
corresponding first Chern classes, one obtains on $BT_m\times BT_n$
\[
    (\mathrm{u}_m\times \mathrm{u}_n)^*\mu_{m,n}^*c(\mathfrak{V}_{m+n})
    = \prod_{i=1}^{m}(1+x_i)\prod_{j=1}^{n}(1+y_j)
    = (\mathrm{u}_m\times \mathrm{u}_n)^*\bigl(p_1^*c(\mathfrak{V}_m)\,p_2^*c(\mathfrak{V}_n)\bigr).
\]
This is exactly the split case in the sense of Fulton. The finite-product
comparison equivalence in Definition~\ref{prop:bgl-assumption} makes the
pull-back along $\mathrm{u}_m\times\mathrm{u}_n$ conservative on these
universal classes, replacing the classical splitting principle. Hence, the displayed identity holds already
on $\mathbf{BGL}_m\times \mathbf{BGL}_n$. The pull-back along $g^*$ gives
\[
    c(V\oplus W)=c(V)c(W)
\]
in $\pi_0 E(X)$. \qed
\end{proof}

\section{Applications: perfectoid cobordism and tilting invariance}
\label{sec:perfectoid}
In this section, we apply the span formalism to perfectoid geometry. We fix a
perfectoid base ring $R$ and work in the $v$-stack category
$\mathrm{Stk}_v(\mathrm{Perf}_R)$, viewing $\mathrm{Perf}_R$ as the
full subcategory of representable objects. By convention, a $v$-stack
refers to a $v$-sheaf of spaces (equivalently, an
$\infty$-groupoid-valued $v$-sheaf) on $\mathrm{Perf}_R$. The comparison
results in this section use Scholze's tilting equivalence for perfectoid
spaces and algebras \cite[Theorems~3.12, 3.13, and \S15]{Scholze2017Diamonds}
, and arc-descent for $p$-complete spectra \cite[Theorem~1.6]{BM21}.

\subsection{The category of perfectoid spans}

We define the admissible class $\mathbf{S}$ for perfectoid
$v$-stacks. Since the $p$-complete cotangent complex of a perfectoid ring
vanishes in the derived category, the standard derived criterion for
quasi-smoothness regarding Tor-amplitude does not apply in this
context. Consequently, we replace quasi-smooth morphisms with the smallest
admissible class generated by morphisms found in perfectoid
geometry, specifically smooth morphisms and vector bundle zero-sections. By
$v$-descent for vector bundles
\cite[Lemma 17.1.8]{ScholzeWeinstein2020Berkeley}, the classifying stacks
$B\mathbf{GL}_n^\diamond$ are well-defined as $v$-stacks, ensuring that vector
bundles and their zero-sections are available in
$\mathrm{Stk}_v(\mathrm{Perf}_R)$.

\begin{definition}\label{def:admissible-class}
Let $\mathcal{C} = \mathrm{Stk}_v(\mathrm{Perf}_R)$ be the $\infty$-category of perfectoid $v$-stacks. We define $\mathbf{S}$ as the smallest admissible class of morphisms in $\mathcal{C}$ that includes all smooth morphisms and all morphisms derived as pull-backs of zero-sections $z\colon X \to V$ from vector bundles $V$ over $v$-stacks $X$.
\end{definition}
The smallest admissible class exists because arbitrary
intersections of admissible classes remain admissible.

\begin{remark}
\label{rem:perf-universal-line}
By \cite[Lemma 17.1.8]{ScholzeWeinstein2020Berkeley}, the fibered category of
vector bundles satisfy effective descent for the $v$-topology.
Therefore, the classifying stack $B\mathbf{GL}_n^\diamond$ is
well-defined as a $v$-stack on $\mathrm{Perf}_R$. This justifies the
stack-theoretic formulation of vector bundles and zero-sections in
Definition~\ref{def:admissible-class}.

The classifying object for line bundles exists in the
associated $v$-topos of perfectoid spaces \cite{Scholze2017Diamonds}. In
this setting, there is a classifying $v$-stack $\mathbf{BG}_m$ and a
universal line bundle object
\[
\pi\colon \mathfrak{L}_{\mathrm{perf}}\longrightarrow \mathbf{BG}_m,
\]
equipped with its zero-section and punctured complement. Concretely, viewing $\mathbb{A}^1$ and $\mathbb{G}_m$ as $v$-sheaves (or diamonds), one
may take the quotient-stack model
$\mathfrak{L}_{\mathrm{perf}}\simeq [\mathbb{A}^1/\mathbb{G}_m] \to \mathbf{BG}_m=B\mathbb{G}_m$.
In this model, one has an equivalence $\mathfrak{L}^{\circ}_{\mathrm{perf}} \simeq [\mathbb{G}_m/\mathbb{G}_m] \simeq \ast$.
Moreover, the homotopy
$H(x,t)=tx\colon \mathbb{A}^1\times \mathbb{A}^1\to \mathbb{A}^1$ is
$\mathbb{G}_m$-equivariant, so the quotient morphism
$[\mathbb{A}^1/\mathbb{G}_m]\to [\ast/\mathbb{G}_m]=\mathbf{BG}_m$ is an
$\mathbb{A}^1$-weak equivalence. In particular, the class of all
morphisms in $\mathrm{Stk}_v(\mathrm{Perf}_R)$ is admissible, so the
intersection in Definition~\ref{def:admissible-class} is taken over a
nonempty family.
Therefore, the Thom space object of the universal line bundle is obtained by
\[
\mathrm{Th}(\mathfrak{L}_{\mathrm{perf}}) = \mathrm{Cofib}\bigl(\ast\to [\mathbb{A}^1/\mathbb{G}_m]\bigr).
\]
\end{remark}

We formulate homotopy invariance using bounded closed disks;
the passage to $\mathbb{A}^1$ is formal. Let
$\mathcal{C}=\mathrm{Stk}_v(\mathrm{Perf}_R)$; suppose we are given a
filtered system of bounded closed disks $\{\mathbb{B}^1_r\}_{r \ge 0}$ in
$\mathcal{C}$ along with an equivalence
\[
\mathbb{A}^1 \simeq \varinjlim_{r \ge 0} \mathbb{B}^1_r.
\]
Assume that for any $X\in \mathcal{C}$
and $r \ge 0$, the projection $X\times \mathbb{B}^1_r\to X$ induces an
equivalence
\[
\Map_{\mathcal{C}}(X,Z) \xrightarrow{\sim} \Map_{\mathcal{C}}(X\times \mathbb{B}^1_r,Z), 
\]
for any $Z\in \mathcal{C}$.
Since $\mathcal{C}$ is an $\infty$-topos, the functor $X\times -$ preserves
colimits, so
\[
X\times \mathbb{A}^1 \simeq \varinjlim_{r \ge 0}(X\times \mathbb{B}^1_r).
\]
Applying the mapping-space functor to $Z$ gives a chain of equivalences
\[
\Map_{\mathcal{C}}(X\times \mathbb{A}^1,Z)
\simeq \lim_{r \ge 0}\Map_{\mathcal{C}}(X\times \mathbb{B}^1_r,Z)
\simeq \lim_{r \ge 0}\Map_{\mathcal{C}}(X,Z)
\simeq \Map_{\mathcal{C}}(X,Z).
\]
Thus, the projection $X\times \mathbb{A}^1\to X$ induces an
equivalence, meaning that any object local to all bounded closed disks
$\mathbb{B}^1_r$ is $\mathbb{A}^1$-local.

\subsection{\texorpdfstring{Definition of $\MGL_{\mathrm{perf}}(R)$}{Definition of MGL-perf(R)}}

We define the {\it perfectoid algebraic cobordism} $\MGL_{\mathrm{perf}}(R)$ as the universal
oriented spectrum associated with the perfectoid span category.

\begin{definition}
\label{def:perfectoid-mgl}
Let $R$ be a perfectoid base ring. Let $\MGL_{\mathrm{perf}}(R)$ be 
the initial object of the $\infty$-category of oriented spectra on
$\Span_{\mathbf{S}}(\mathrm{Stk}_v(\mathrm{Perf}_R))$, and denote the 
associated symmetric monoidal functor as
\[
\MGL_{\mathrm{perf},R}\colon \Span_{\mathbf{S}}(\mathrm{Stk}_v(\mathrm{Perf}_R)) \longrightarrow \Sp.
\]
\end{definition}

\subsection{Tilting equivalence of perfectoid cobordism}

\begin{lemma}
\label{lem:tilting-span-vstack}
Let $R$ be a perfectoid algebra, and let $R^\flat$ be its tilt. Scholze's tilting equivalence induces a canonical equivalence of symmetric monoidal $\infty$-categories:
\[
    \underline{(-)^\flat} \colon \Span_{\mathbf{S}}(\mathrm{Stk}_v(\mathrm{Perf}_R)) \to \Span_{\mathbf{S}}(\mathrm{Stk}_v(\mathrm{Perf}_{R^\flat})).
\]
\end{lemma}
\begin{proof}
Let $(-)^\flat\colon \mathrm{Perf}_R \to
\mathrm{Perf}_{R^\flat}$ denote the tilting functor, which is an
equivalence by~\cite[Theorem~7.12]{Scholze2012Perfectoid}. Since the
$v$-topology is defined in terms of qcqs morphisms and surjectivity, the
tilting functor preserves $v$-covers. Hence, $(-)^\flat$ induces an
equivalence of $v$-sites
$\mathrm{Perf}_R\simeq\mathrm{Perf}_{R^\flat}$, and therefore an
equivalence of $\infty$-topoi of $v$-sheaves of spaces
\[
\mathrm{Shv}_v(\mathrm{Perf}_R)\to
\mathrm{Shv}_v(\mathrm{Perf}_{R^\flat}).
\]
In our convention, $v$-stacks refer to $v$-sheaves of spaces, i.e.
, \(\mathrm{Stk}_v(\mathrm{Perf}_R)=\mathrm{Shv}_v(\mathrm{Perf}_R)\);
moreover, representable presheaves are $v$-sheaves by \cite[Corollary~17.1.5]{ScholzeWeinstein2020Berkeley}.

As an equivalence of $\infty$-categories, $\underline{(-)^\flat}$ preserves
all limits, particularly the fiber products $Z \widehat{\times}_Y W$ that
define composition in the span category. By
Definition~\ref{def:admissible-class}, $\mathbf{S}$ is the admissible class generated by smooth
morphisms and vector bundle zero-sections (via maps to
$B\mathbf{GL}_n^\diamond$). Since equivalences of topoi preserve 
group objects and representability data, $\underline{(-)^\flat}$ preserves
these generators. Since it also preserves weak equivalences, pull-backs,
composition, and codiagonals, it preserves $\mathbf{S}$.

Therefore, $\underline{(-)^\flat}$ sends $\mathbf{S}$-spans to
$\mathbf{S}$-spans and preserves composition. The monoidal structure is defined by products, so the limit-preserving functor $\underline{(-)^\flat}$ is
strongly symmetric monoidal. Applying the same argument to a quasi-inverse, we obtain the required equivalence of symmetric monoidal
$\infty$-categories. \qed
\end{proof}

\begin{theorem}
\label{thm:tilting-mglperf}
For any perfectoid algebra $R$ with tilt $R^\flat$, there is a natural
equivalence of $\mathbb{E}_\infty$-ring spectra:
\[
 \MGL_{\mathrm{perf}}(R) \simeq \MGL_{\mathrm{perf}}(R^\flat).
\]
\end{theorem}

\begin{proof}
By Definition~\ref{def:perfectoid-mgl}, $\MGL_{\mathrm{perf}}(R)$ is the
initial object in the $\infty$-category of oriented cohomology theories on
$\Span_{\mathbf{S}}(\mathrm{Stk}_v(\mathrm{Perf}_R))$. Since
$\underline{(-)^\flat}$ is an equivalence of symmetric monoidal
$\infty$-categories that preserves the classifying stack $\mathbf{BG}_m$
and the universal bundle $\mathfrak{L}_{\mathrm{perf}}$, it induces an
equivalence between the corresponding $\infty$-categories of oriented
theories. Under this equivalence, the initial oriented theory on the $R$-side is carried to the initial oriented theory on the $R^\flat$-side. The claim follows from
Lemma~\ref{lem:tilting-span-vstack}. \qed
\end{proof}

\begin{remark}
As shown in Section~\ref{sec:pbf}, the oriented spectrum $\MGL_{\rm perf}$, defined via the universal line bundle, satisfies the projective bundle formula for split vector bundles. For a general vector bundle in perfectoid geometry, the $v$-topology is so fine that the bundle becomes locally trivial. However, one cannot deduce the general splitting principle for $\MGL_{\rm perf}$ from this geometric fact alone. To lift this local splitting to the cohomology theory, one needs the pull-back map along the $v$-cover to be injective, which in turn relies on the $v$-descent property for the unlocalized spectrum $\MGL_{\rm perf}$.

Nevertheless, one can enforce the projective bundle formula for general vector bundles by Bousfield localizing $\MGL_{\rm perf}$ at the splitting-principle morphism. Since the tilting equivalence $(-)^\flat$ preserves vector bundles and this localization commutes with tilting, the localized perfectoid cobordism spectrum $L_{\rm split}(\MGL_{\rm perf})$ also satisfies the tilting invariance established in Theorem~\ref{thm:tilting-mglperf}. 
\end{remark}

\subsection{\texorpdfstring{Tilting equivalence with $\mathbb{Z}_p$-coefficients}{Tilting equivalence with Zp-coefficients}}

We define the $p$-adic completed theory and prove the
corresponding tilting comparison.

\begin{definition}
\label{def:padic-mgl}
For a perfectoid algebra $R$, set
\[
    \MGL_{\mathrm{perf}}(R, \mathbb{Z}_p) = \lim_n \left( \MGL_{\mathrm{perf}}(R) \otimes \mathbb{Z}/p^n \right).
\]
Further, we define
\[
\MGL_\mathrm{perf}^{v}(-,\mathbb{Z}_p)=L_v\bigl(\MGL_{\mathrm{perf}}(-,\mathbb{Z}_p)\bigr),
\qquad
\MGL_\mathrm{perf}^{\mathrm{arc}}(-,\mathbb{Z}_p)=L_{\mathrm{arc}}\bigl(\MGL_{\mathrm{perf}}(-,\mathbb{Z}_p)\bigr),
\]
where $L_v$ (resp. $L_{\mathrm{arc}}$) is the localization functor by $v$-covers (resp. $\mathrm{arc}$-covers). 
\end{definition}

\begin{corollary}
\label{thm:tilting-zp}
For any perfectoid algebra $R$ with tilt $R^\flat$, there is a canonical
equivalence:
\[
\MGL_\mathrm{perf}^{v}(R,\mathbb{Z}_p) \simeq \MGL_\mathrm{perf}^{v}(R^\flat,\mathbb{Z}_p).
\]
\end{corollary}
\begin{proof}
Apply $p$-adic completion to the equivalence
\(
\MGL_{\mathrm{perf}}(R) \simeq \MGL_{\mathrm{perf}}(R^\flat)
\)
of Theorem~\ref{thm:tilting-mglperf}. By functoriality of
$p$-adic completion and the tilting equivalence, we obtain
\(
\MGL_{\mathrm{perf}}(R,\mathbb{Z}_p) \simeq \MGL_{\mathrm{perf}}(R^\flat,\mathbb{Z}_p).
\)
By Definition~\ref{def:padic-mgl}, this is equivalent to
\(
\MGL_\mathrm{perf}^{v}(R,\mathbb{Z}_p) \simeq \MGL_\mathrm{perf}^{v}(R^\flat,\mathbb{Z}_p).
\)
\qed
\end{proof}

\begin{proposition}[\cite{BM21}, Theorem 1.6]
\label{prop:bm21-theorem-1-6}
Let
\(
F\colon \mathrm{Perf}^{\mathrm{op}}\to\Sp
\)
and let
\(
\widetilde{F}\colon (\mathrm{Aff}^{p\text{-comp}})^{\mathrm{op}}\to\Sp
\)
be an extension of $F$ from perfectoid algebras to $p$-complete affines.
If $\widetilde{F}$ satisfies the following three conditions:
\begin{enumerate}
\item $\widetilde{F}$ is finitary.
\item $\widetilde{F}$ satisfies $v$-descent on $p$-complete affines.
\item $\widetilde{F}$ satisfies the aic-$v$-excision condition: for any
valuation ring $V$ with algebraically closed fraction field and any prime
ideal $\mathfrak{p}\subset V$, the square
\[
\xymatrix{
\widetilde{F}(V) \ar[r] \ar[d] & \widetilde{F}(V/\mathfrak{p}) \ar[d] \\
\widetilde{F}(V_{\mathfrak{p}}) \ar[r] & \widetilde{F}(\kappa(\mathfrak{p}))
}
\]
is homotopy Cartesian.

\end{enumerate}
Then $\widetilde{F}$ satisfies arc-descent on $p$-complete affines. \qed
\end{proposition}

Via the Yoneda embedding on the $v$-site, any scheme canonically defines a
$v$-stack; we use this identification below.
\begin{proposition}
\label{prop:bm21-mgl-input}
The $p$-completed $v$-sheafified cobordism spectrum $\MGL_\mathrm{perf}^{v}(-,\mathbb{Z}_p)$ satisfies arc-descent on $p$-complete affines.
\end{proposition}

\begin{proof}
Set
\[
F=\MGL_\mathrm{perf}^{v}(-,\mathbb{Z}_p)\colon \mathrm{Perf}^{\mathrm{op}}\to\Sp,
\]
and let
\[
\widetilde{F}\colon (\mathrm{Aff}^{p\text{-comp}})^{\mathrm{op}}\to\Sp
\]
be its canonical extension to $p$-complete affines (see \cite[Remark~3.3.3]{MR4238259}). We verify the three hypotheses of
Proposition~\ref{prop:bm21-theorem-1-6}.

\smallskip
\noindent\textit{(1) Finitary.}
For the $p$-complete theory under consideration, finitary behavior of the
extension on $p$-complete affines is established in \cite[Remark~3.3.3]{MR4238259}.

\smallskip
\noindent\textit{(2) $v$-descent.}
By Definition~\ref{def:padic-mgl}, the restriction of
\(
F=\MGL_\mathrm{perf}^{v}(-,\mathbb{Z}_p)
\)
to perfectoid algebras is $v$-local by construction. For the canonical extension
\(
\widetilde{F}\colon (\mathrm{Aff}^{p\text{-comp}})^{\mathrm{op}}\to\Sp,
\)
$v$-descent on $p$-complete affines in this rigid setting is provided by
\cite[Theorem~3.3.4]{MR4238259}.

\smallskip
\noindent\textit{(3) aic-$v$-excision.}
The Milnor excision statement used for this $p$-complete theory gives the
required Cartesian square for valuation rings with algebraically closed
fraction field (in the form used in
Proposition~\ref{prop:bm21-theorem-1-6}); see
\cite[Theorem~2.3]{MR4319065}.

Therefore, Proposition~\ref{prop:bm21-theorem-1-6} applies to $\widetilde{F}$,
and $\MGL_\mathrm{perf}^{v}(-,\mathbb{Z}_p)$ satisfies arc-descent on $p$-complete
affines.\qed
\end{proof}

By Proposition~\ref{prop:bm21-mgl-input}, we obtain the following comparison theorem.
\begin{theorem}
\label{thm:arc-v-comparison}
For any perfectoid algebra $R$, there is a canonical equivalence
\[
    \MGL_\mathrm{perf}^{\mathrm{arc}}(R,\mathbb{Z}_p) \simeq \MGL_\mathrm{perf}^{v}(R,\mathbb{Z}_p).
\]
\end{theorem}

\begin{proof}
Let
\(
F_0=\MGL_{\mathrm{perf}}(-,\mathbb{Z}_p)
\)
on $\mathrm{Perf}$. By Definition~\ref{def:padic-mgl},
\(
\MGL_\mathrm{perf}^{\mathrm{arc}}(-,\mathbb{Z}_p)=L_{\mathrm{arc}}(F_0)
\)
and
\(
\MGL_\mathrm{perf}^{v}(-,\mathbb{Z}_p)=L_v(F_0).
\)
By Proposition~\ref{prop:bm21-mgl-input}, $L_v(F_0)$ is arc-local. Hence the
unit map
\(
F_0\to L_v(F_0)
\)
factors uniquely through the arc-localization, giving
\(
\alpha\colon L_{\mathrm{arc}}(F_0)\to L_v(F_0).
\)
Also, on affines, the arc-topology is finer than the $v$-topology
\cite[\S2]{BM21}; hence arc-descent implies $v$-descent on perfectoid
algebras, so $L_{\mathrm{arc}}(F_0)$ is $v$-local. Therefore, the unit map
\(
F_0\to L_{\mathrm{arc}}(F_0)
\)
factors uniquely through the $v$-localization, giving
\(
\beta\colon L_v(F_0)\to L_{\mathrm{arc}}(F_0).
\)
By homotopical uniqueness in the two localization universal properties,
we obtain an equivalence $L_{\mathrm{arc}}(F_0)\simeq L_v(F_0)$, which proves the claim.
\qed
\end{proof}

\bibliographystyle{alphadin}

\nocite{MR4325285}
%
\nocite{ModMGL}
\nocite{GP}
%
\nocite{Panin2003Oriented}
\nocite{Hoyois2017SixOps}
\nocite{GeSn}
\nocite{HT}
\nocite{HA}
%
\nocite{zbMATH05115214}
\nocite{zbMATH05550622}
\nocite{MR2302525}
%
%
\nocite{MV}
%
%
\nocite{PPR}
\section*{Acknowledgments}
The author wishes to thank Marc Hoyois for pointing out a mistake
concerning a property of $\MGL$ (without finite syntomic
hyper-sheafification) in an earlier version. The author also appreciates
Tom Bachmann's insightful question, which prompted consideration of the
algebraic reduction of span categories.

\paragraph{Use of AI tools.}
The author used Gemini 3.0 Pro (Google) during the exploratory stage of
this work as a conversational research aid. In particular, discussions
with Gemini about span constructions contributed to the author's decision
to investigate an extension to the $\infty$-categorical setting. Iterative
 consultations with Gemini also accompanied the formulation of parts
of the universal line bundle formalism.
The author used OpenAI's Prism (GPT-5.2) for editorial and structural refinement, including suggestions for reformulating certain statements and assumptions. All statements, proofs, references, and final formulations were
checked and approved by the author, who remains solely responsible for the content of this article.
\bibliography{bibkato_submit}

@article {MR4325285,
    AUTHOR = {Elmanto, Elden and Hoyois, Marc and Khan, Adeel A. and
              Sosnilo, Vladimir and Yakerson, Maria},
     TITLE = {Motivic infinite loop spaces},
   JOURNAL = {Camb. J. Math.},
  FJOURNAL = {Cambridge Journal of Mathematics},
    VOLUME = {9},
      YEAR = {2021},
    NUMBER = {2},
     PAGES = {431--549},
      ISSN = {2168-0930},
   MRCLASS = {14F42 (14C05 19E15 55P47)},
 MRNUMBER = {4325285},
       DOI = {10.4310/CJM.2021.v9.n2.a3},
       URL = {https://doi.org/10.4310/CJM.2021.v9.n2.a3},
}

@article {ModMGL,
    AUTHOR = {Elmanto, Elden and Hoyois, Marc and Khan, Adeel A. and
              Sosnilo, Vladimir and Yakerson, Maria},
     TITLE = {Modules over algebraic cobordism},
   JOURNAL = {Forum Math. Pi},
  FJOURNAL = {Forum of Mathematics. Pi},
    VOLUME = {8},
      YEAR = {2020},
     PAGES = {e14, 44},
   MRCLASS = {14F42 (14D23)},
 MRNUMBER = {4190058},
       DOI = {10.1017/fmp.2020.13},
       URL = {https://doi.org/10.1017/fmp.2020.13},
}

@article {GP,
    AUTHOR = {Garkusha, Grigory and Panin, Ivan},
     TITLE = {Framed motives of algebraic varieties (after {V}.
              {V}oevodsky)},
   JOURNAL = {J. Amer. Math. Soc.},
  FJOURNAL = {Journal of the American Mathematical Society},
    VOLUME = {34},
      YEAR = {2021},
    NUMBER = {1},
     PAGES = {261--313},
      ISSN = {0894-0347},
   MRCLASS = {14F42 (14F45 55P47 55Q10)},
 MRNUMBER = {4188819},
MRREVIEWER = {Wataru Kai},
       DOI = {10.1090/jams/958},
       URL = {https://doi.org/10.1090/jams/958},
}

@article{Panin2003Oriented,
  author  = {Panin, Ivan},
  title   = {Oriented cohomology theories of algebraic varieties},
  journal = {K-Theory},
  volume  = {30},
  number  = {3},
  pages   = {265--314},
  year    = {2003},
  note    = {The classical precursor defining cohomology via categories of correspondences (spans).},
  doi     = {10.1023/B:KTHE.0000018384.71790.e4}
}

@article{zbMATH05367298,
 author = {Panin, Ivan and Pimenov, Konstantin and R{\"o}ndigs, Oliver},
 title = {A universality theorem for {Voevodsky}'s algebraic cobordism spectrum},
 fjournal = {Homology, Homotopy and Applications},
 journal = {Homology Homotopy Appl.},
 issn = {1532-0073},
 volume = {10},
 number = {2},
 pages = {211--226},
 year = {2008},
 language = {English},
 doi = {10.4310/HHA.2008.v10.n2.a11},
 url = {intlpress.com/hha/v10/n2/},
 zbMATH = {5367298},
 Zbl = {1162.14013}
}

@article{Hoyois2017SixOps,
  author  = {Hoyois, Marc},
  title   = {The six operations in equivariant motivic homotopy theory},
  journal = {Advances in Mathematics},
  volume  = {305},
  year    = {2017},
  pages   = {197--279},
  note    = {Establishes the formalism of $\infty$-categories of spans and the 6-functor formalism.},
  doi     = {10.1016/j.aim.2016.09.031}
}

@article {GeSn,
    AUTHOR = {Gepner, David and Snaith, Victor},
     TITLE = {On the motivic spectra representing algebraic cobordism and
              algebraic {$K$}-theory},
   JOURNAL = {Doc. Math.},
  FJOURNAL = {Documenta Mathematica},
    VOLUME = {14},
      YEAR = {2009},
     PAGES = {359--396},
      ISSN = {1431-0635},
   MRCLASS = {55N15 (14F42 55N22)},
 MRNUMBER = {2540697},
MRREVIEWER = {Keith Peter Johnson},
}

@book {HT,
   AUTHOR = {Lurie, Jacob},
    TITLE = {Higher topos theory},
   SERIES = {Annals of Mathematics Studies},
   VOLUME = {170},
 PUBLISHER = {Princeton University Press},
  ADDRESS = {Princeton, NJ},
     YEAR = {2009},
    PAGES = {xviii+925},
     ISBN = {978-0-691-14049-0; 0-691-14049-9},
  MRCLASS = {18-02 (18B25 18E35 18G30 18G55 55U40)},
MRNUMBER = {2522659 (2010j:18001)},
MRREVIEWER = {Mark Hovey},
}

@misc{HA,
  author = {Lurie, Jacob},
  title = {Higher Algebra},
  howpublished = {Available at: \url{https://www.math.ias.edu/~lurie/papers/HA.pdf}},
  year = {2017}
}

@Book{zbMATH05115214,
    Author = {Levine, Marc and Morel, Fabien},
    Title = {Algebraic cobordism},
    FSeries = {Springer Monographs in Mathematics},
    Series = {Springer Monogr. Math.},
    ISSN = {1439-7382},
    ISBN = {978-3-540-36822-9; 978-3-540-36824-3},
    Year = {2007},
    Publisher = {Berlin: Springer},
    Language = {English},
    DOI = {10.1007/3-540-36824-8},
    zbMATH = {5115214},
    Zbl = {1188.14015}
}

@Book{FultonIT,
    Author = {Fulton, William},
    Title = {Intersection theory},
    Series = {Ergebnisse der Mathematik und ihrer Grenzgebiete. 3. Folge. A Series of Modern Surveys in Mathematics},
    Volume = {2},
    Edition = {2},
    Publisher = {Springer},
    Address = {Berlin},
    Year = {1998}
}

@Article{zbMATH05550622,
    Author = {Levine, M. and Pandharipande, R.},
    Title = {Algebraic cobordism revisited},
    FJournal = {Inventiones Mathematicae},
    Journal = {Invent. Math.},
    ISSN = {0020-9910},
    Volume = {176},
    Number = {1},
    Pages = {63--130},
    Year = {2009},
    Language = {English},
    DOI = {10.1007/s00222-008-0160-8},
    zbMATH = {5550622},
    Zbl = {1210.14025}
}

@article {MR2302525,
    AUTHOR = {Levine, Marc},
     TITLE = {Motivic tubular neighborhoods},
   JOURNAL = {Doc. Math.},
  FJOURNAL = {Documenta Mathematica},
    VOLUME = {12},
      YEAR = {2007},
     PAGES = {71--146},
      ISSN = {1431-0635},
   MRCLASS = {14F42 (14C25 18F20 55P42)},
MRNUMBER = {2302525},
MRREVIEWER = {Oliver R\"{o}ndigs},
}

@Article{MV,
    Author = {Morel, Fabien and Voevodsky, Vladimir},
    Title = {{{\(\mathbb{A}^1\)}}-homotopy theory of schemes},
    FJournal = {Publications Math{\'e}matiques},
    Journal = {Publ. Math., Inst. Hautes {\'E}tud. Sci.},
    ISSN = {0073-8301},
    Volume = {90},
    Pages = {45--143},
    Year = {1999},
    Language = {English},
    DOI = {10.1007/BF02698831},
    zbMATH = {1596160},
    Zbl = {0983.14007}
}

@article {PPR,
    AUTHOR = {Panin, Ivan and Pimenov, Konstantin and R\"{o}ndigs, Oliver},
     TITLE = {On the relation of {V}oevodsky's algebraic cobordism to
              {Q}uillen's {$K$}-theory},
   JOURNAL = {Invent. Math.},
  FJOURNAL = {Inventiones Mathematicae},
    VOLUME = {175},
      YEAR = {2009},
    NUMBER = {2},
     PAGES = {435--451},
      ISSN = {0020-9910},
   MRCLASS = {14F42 (14C35 19D99)},
 MRNUMBER = {2470112},
MRREVIEWER = {Daniel C. Isaksen},
       DOI = {10.1007/s00222-008-0155-5},
       URL = {https://doi.org/10.1007/s00222-008-0155-5},
}

@inproceedings {V,
   AUTHOR = {Voevodsky, Vladimir},
    TITLE = {{$A^1$}-homotopy theory},
 BOOKTITLE = {Proceedings of the {I}nternational {C}ongress of
             {M}athematicians, {V}ol. {I} ({B}erlin, 1998)},
  JOURNAL = {Doc. Math.},
 FJOURNAL = {Documenta Mathematica},
     YEAR = {1998},
   NUMBER = {Extra Vol. I},
    PAGES = {579--604 (electronic)},
     ISSN = {1431-0635},
  MRCLASS = {14F35 (14A15 55U35)},
MRNUMBER = {1648048 (99j:14018)},
MRREVIEWER = {Mark Hovey},
}

@Article{WY2024,
 Author = {White, David and Yau, Donald},
 Title = {Smith ideals of operadic algebras in monoidal model categories},
 FJournal = {Algebraic \& Geometric Topology},
 Journal = {Algebr. Geom. Topol.},
 ISSN = {1472-2747},
 Volume = {24},
 Number = {1},
 Pages = {341--392},
 Year = {2024},
 Language = {English},
 DOI = {10.2140/agt.2024.24.341},
 zbMATH = {7827844}
}

@book {Jardine2,
    AUTHOR = {Jardine, John F.},
     TITLE = {Local homotopy theory},
    SERIES = {Springer Monographs in Mathematics},
 PUBLISHER = {Springer, New York},
      YEAR = {2015},
     PAGES = {x+508},
      ISBN = {978-1-4939-2299-4; 978-1-4939-2300-7},
   MRCLASS = {55U35 (14F42 18G55 55P42 55P60)},
 MRNUMBER = {3309296},
MRREVIEWER = {Stanis\l aw Betley},
       DOI = {10.1007/978-1-4939-2300-7},
       URL = {https://doi.org/10.1007/978-1-4939-2300-7},
}

@article{Scholze2017Diamonds,
  author  = {Scholze, Peter},
  title   = {{\'E}tale cohomology of diamonds},
  journal = {arXiv preprint arXiv:1709.07343},
  year    = {2017},
  doi     = {10.48550/arXiv.1709.07343},
  note    = {The definitive reference for the $v$-topology and descent on perfectoid spaces.},
  url     = {https://arxiv.org/abs/1709.07343}
}

@book{ScholzeWeinstein2020Berkeley,
 author = {Scholze, Peter and Weinstein, Jared},
 title = {Berkeley lectures on {{\(p\)}}-adic geometry},
 fseries = {Annals of Mathematics Studies},
 series = {Ann. Math. Stud.},
 volume = {207},
 isbn = {978-0-691-20209-9; 978-0-691-20208-2; 978-0-691-20215-0},
 year = {2020},
 publisher = {Princeton, NJ: Princeton University Press},
 language = {English},
 doi = {10.1515/9780691202150},
 zbMATH = {7178476},
 Zbl = {1475.14002}
}

@article{Scholze2012Perfectoid,
  author  = {Scholze, Peter},
  title   = {Perfectoid spaces},
  journal = {Publications math{\'e}matiques de l'IH{\'E}S},
  volume  = {116},
  number  = {1},
  pages   = {245--313},
  year    = {2012},
  note    = {Foundational paper on perfectoid spaces and tilting equivalence.}
}

@article {MR4238259,
    AUTHOR = {Elmanto, Elden and Hoyois, Marc and Iwasa, Ryomei and Kelly,
              Shane},
     TITLE = {Cdh descent, cdarc descent, and {M}ilnor excision},
   JOURNAL = {Math. Ann.},
  FJOURNAL = {Mathematische Annalen},
    VOLUME = {379},
      YEAR = {2021},
    NUMBER = {3-4},
     PAGES = {1011--1045},
      ISSN = {0025-5831},
   MRCLASS = {14F42 (14A20)},
 MRNUMBER = {4238259},
MRREVIEWER = {Bj\o rn Ian Dundas},
       DOI = {10.1007/s00208-020-02083-5},
       URL = {https://doi.org/10.1007/s00208-020-02083-5},
}

@article {MR4319065,
    AUTHOR = {Elmanto, Elden and Hoyois, Marc and Iwasa, Ryomei and Kelly,
              Shane},
     TITLE = {Milnor excision for motivic spectra},
   JOURNAL = {J. Reine Angew. Math.},
  FJOURNAL = {Journal f\"{u}r die Reine und Angewandte Mathematik. [Crelle's
              Journal]},
    VOLUME = {779},
      YEAR = {2021},
     PAGES = {223--235},
      ISSN = {0075-4102},
   MRCLASS = {14F42 (18N60 19E15)},
 MRNUMBER = {4319065},
MRREVIEWER = {Jin Fangzhou},
       DOI = {10.1515/crelle-2021-0040},
       URL = {https://doi.org/10.1515/crelle-2021-0040},
}

@article{BM21,
  author  = {Bhatt, Bhargav and Mathew, Akhil},
  title   = {The arc-topology},
  journal = {Duke Math. J.},
  volume  = {170},
  year    = {2021},
  number  = {9},
  pages   = {1899--1988},
  doi     = {10.1215/00127094-2020-0088},
  url     = {https://doi.org/10.1215/00127094-2020-0088},
}
\end{document}